\documentclass[11pt]{article}%
\usepackage{amssymb}
\usepackage{amsfonts}
\usepackage{amsmath}
\usepackage{graphicx}%
\providecommand{\U}[1]{\protect\rule{.1in}{.1in}}
\newtheorem{theorem}{Theorem}

\newtheorem{corollary}{Corollary}

\newcounter{article}

\begin{document}

\title{Graph eigenvectors, fundamental weights and centrality metrics for nodes in networks}
\author{P. Van Mieghem\thanks{ Faculty of Electrical Engineering, Mathematics and
Computer Science, P.O Box 5031, 2600 GA Delft, The Netherlands; \emph{email}:
P.F.A.VanMieghem@tudelft.nl }}
\date{Delft University of Technology\\
14 February 2014 (v1)\\
10 December 2014 (v2)\\
8 August 2015 (v3)\\
15 March 2016 (v4)}
\maketitle

\begin{abstract}
Several expressions for the $j$-th component $\left(  x_{k}\right)  _{j}$ of
the $k$-th eigenvector $x_{k}$ of a symmetric matrix $A$ belonging to
eigenvalue $\lambda_{k}$ and normalized as $x_{k}^{T}x_{k}=1$ are presented.
In particular, the expression
\[
\left(  x_{k}\right)  _{j}^{2}=-\frac{1}{c_{A}^{\prime}\left(  \lambda
_{k}\right)  }\det\left(  A_{\backslash\left\{  j\right\}  }-\lambda
_{k}I\right)
\]
where $c_{A}\left(  \lambda\right)  =\det\left(  A-\lambda I\right)  $ is the
characteristic polynomial of $A$, $c_{A}^{\prime}\left(  \lambda\right)
=\frac{dc_{A}\left(  \lambda\right)  }{d\lambda}$ and $A_{\backslash\left\{
j\right\}  }$ is obtained from $A$ by removal of row $j$ and column $j$,
suggests us to consider the square eigenvector component as a graph centrality
metric for node $j$ that reflects the impact of the removal of node $j$ from
the graph at an eigenfrequency/eigenvalue $\lambda_{k}$ of a graph related
matrix (such as the adjacency or Laplacian matrix). Removal of nodes in a
graph relates to the robustness of a graph. The set of such nodal centrality
metrics, the squared eigenvector components $\left(  x_{k}\right)  _{j}^{2}$
of the adjacency matrix over all eigenvalue $\lambda_{k}$ for each node $j$,
is \textquotedblleft ideal\textquotedblright\ in the sense of being complete,
\emph{almost} uncorrelated and mathematically precisely defined and
computable. Fundamental weights (column sum of $X$) and dual fundamental
weights (row sum of $X$) are introduced as spectral metrics that condense
information embedded in the orthogonal eigenvector matrix $X$, with elements
$X_{ij}=\left(  x_{j}\right)  _{i}$.

In addition to the criterion \textquotedblleft If the algebraic connectivity
is positive, then the graph is connected\textquotedblright, we found an
alternative condition: \textquotedblleft If $\min_{1\leq k\leq N}\left(
\lambda_{k}^{2}(A)\right)  =d_{\min}$, then the graph is
disconnected.\textquotedblright

\end{abstract}

\section{Introduction}

\label{sec_introduction}Generally, nodal centrality metrics quantify the
\textquotedblleft importance\textquotedblright\ of a node\footnote{The
importance of a link in $G$ can be assessed as the importance of a node in the
corresponding line graph $l\left(  G\right)  $, defined in \cite[p.
17-21]{PVM_graphspectra}.} in a network or how \textquotedblleft
central\textquotedblright\ a node is in the graph. Many quantifiers of nodal
\textquotedblleft importance\textquotedblright\ have been proposed, that are
reviewed in
\cite{Boccaletti_PhysicsReports2006,PhD_Javier_MartinHernandez,Estrada2012}.
Perhaps, the simplest -- both in meaning as well as in computation -- is the
degree of a node defined as the number of direct neighbors of a node in the
network. Relevant questions such as \textquotedblleft What is the most
influential node in a social networks?\textquotedblright%
\ \cite{Kitsak_Stanley_nature2010}\ and \textquotedblleft What is the most
vulnerable node when attacked or removed?\textquotedblright\ are difficult to
answer, because a precise translation of \textquotedblleft
influence\textquotedblright\ or \textquotedblleft
vulnerability\textquotedblright\ in terms of computable quantities, called
metrics \cite{PVM_PAComplexNetsCUP}, of the graph is needed. Nodal
\textquotedblleft importance\textquotedblright\ often depends on the process
on the network, which then further specifies the precise meaning of importance
with respect to that process. For example, in epidemics on networks
\cite{PVM_RMP_epidemics2014}, nodal importance (here vulnerability) can be
defined as the long-run probability that a node is infected
\cite{PVM_N_intertwined_Computing2011}, given an effective infection rate
$\tau$ of the virus. The most \textquotedblleft influential\textquotedblright%
\ spreader can be defined as the fastest spreader, that, when initially
injected with information, reaches in the shortest time the metastable
fraction of infected nodes, again given an effective infection rate $\tau$.
Both the nodal ranking in vulnerability and the fastest spreader change with
effective infection rate $\tau$, clearly illustrating that only topological
metrics are inadequate to determine the \textquotedblleft most
important\textquotedblright\ node.  

Besides the precise definition, meaning and applicability or usefulness of a
graph metric, a number of other issues appear as elaborated in
\cite{PVM_R_model}: How many metrics are needed to compare graphs? How
strongly is a set of two metrics correlated? How difficult is the computation
of the metric and how much information of the network is required (only local
information as the degree or global information as for the diameter)? In most
cases, more than one metric is needed to quantify the desired
\textquotedblleft importance\textquotedblright. For example, a high-degree
node of which all neighbors have degree 1 and one neighbor has degree 2, is
vulnerable to be disconnected from the remainder of the network, in spite of
its high degree. When multiple metrics are chosen, they should be as
independent or orthogonal as possible, because strongly correlated metrics can
be combined to a single one, since they all reflect the same type of
\textquotedblleft importance\textquotedblright\ as illustrated in
\cite{PVM_correlation_metrics_JSTAT2011}.

Here, we take a different view. We present a complete set of orthogonal
centrality metrics and try to interpret what type of properties in the network
they may characterize or quantify. As reviewed in Appendix
\ref{sec_introduction_eigenvectors_eigenvalues}, a non-zero vector $x\left(
\lambda\right)  $ only satisfies the eigenvalue equation%
\[
Ax\left(  \lambda\right)  =\lambda\ x\left(  \lambda\right)
\]
if the real number $\lambda$, which we can interpret as a \textquotedblleft
frequency\textquotedblright, is an eigenvalue of $A$ such that $x_{k}=x\left(
\lambda_{k}\right)  $ is the eigenvector at eigenfrequency $\lambda
=\lambda_{k}$. We normalize $x_{k}$ so that $x_{k}^{T}x_{k}=1$, according to
the \emph{first} \cite{PVM_double_orthogonality} orthogonality equations
(\ref{first_orthogonality_equations}) and denote the $j$-th eigenvector
component by $\left(  x_{k}\right)  _{j}$, where the index $j$ refers to nodes
and the index $k$ to eigenfrequencies. Three different expressions
(\ref{squared_eigenvector_component_xk_j_charc_pol}),
(\ref{i_th_component_xk_square}), and (\ref{(x_m)^2_char_pol}) for the
\emph{square} of the $j$-th component of the $k$-th eigenvector $\left(
x_{k}\right)  _{j}$ of the adjacency matrix $A$ belonging to eigenvalue
$\lambda_{k}$ are presented. The determinantal expression
(\ref{squared_eigenvector_component_xk_j_charc_pol}) is derived in Section
\ref{sec_determinantal_square_eigenvector_component}, essentially using merely
linear algebra. Section \ref{sec_determinantal_square_eigenvector_component}
further interprets expression
(\ref{squared_eigenvector_component_xk_j_charc_pol}) for $\left(
x_{k}\right)  _{j}^{2}$ as the \emph{impact of the removal of node }$j$\emph{
from }$G$\emph{ at eigenfrequency }$\lambda_{k}$ of a symmetric graph matrix
(such as the adjacency matrix or the Laplacian). Strongly based on the
eigenvalue equation of the adjacency matrix $A$, Section
\ref{sec_squared_eigenvalue_equation} derives the second expression
(\ref{i_th_component_xk_square}) for $\left(  x_{k}\right)  _{j}^{2}$. Several
bounds are given, of which some extend earlier published bounds. The third,
walk-based expression (\ref{(x_m)^2_char_pol}) for $\left(  x_{k}\right)
_{j}^{2}$ is deferred to Appendix \ref{sec_walk_based_(xk)_j}: for reasons of
completeness, we have incorporated (\ref{(x_m)^2_char_pol}). The elegance of
(\ref{squared_eigenvector_component_xk_j_charc_pol}) illustrates that the
\emph{square} $\left(  x_{k}\right)  _{j}^{2}$ is likely more suited than
$\left(  x_{k}\right)  _{j}$ to explain the behavior of the eigenstructure,
which reminds us to the basic interpretation of quantum mechanics (see e.g.
\cite{Dirac,CohenTannoudji}), where the wave function can be complex, while
its modulus is interpreted as a probability. Unfortunately, as shown in
Section \ref{sec_centrality_vectors_are_dependent} for the adjacency matrix
$A$, the vector $c_{j}=$ $\left(  \left(  x_{1}\right)  _{j}^{2},\left(
x_{2}\right)  _{j}^{2},\ldots,\left(  x_{N}\right)  _{j}^{2}\right)  $ of the
adjacency centrality metrics $\left(  x_{k}\right)  _{j}^{2}$ at all
eigenfrequencies $k$ for node $j$ is not independent (or orthogonal) to
$c_{l}$ for node $l$, which implies that the set of adjacency eigenvector
centrality metrics $\left\{  c_{j}\right\}  _{1\leq j\leq N}$ is not complete!

Section \ref{sec_fundamental_weights} introduces the definitions and basic
properties of the fundamental weights and the dual fundamental weights of a
graph. Fundamental weights and their dual are proposed as possible
condensations of the $N\times N$ orthogonal matrix $X$ containing all
eigenvectors $\left\{  x_{k}\right\}  _{1\leq k\leq N}$ in its columns. The
aim to find a more economical way (i.e. less than $N^{2}$ elements) for $X$,
while not loosing information (i.e. able to reconstruct $X$), started already
with Cvetkovic \cite{Cvetkovic_eigenspaces}, who introduced \textquotedblleft
graph angles\textquotedblright. For a sufficiently large graph, Van Dam and
Haemers \cite{VanDam_Haemers2003} have argued that the set of all eigenvalues
alone (thus ignoring eigenvectors or $X$) is a unique fingerprint or signature
of the graph. For exact graph reconstruction and storage of networks, the most
condensed form (i.e. least number of bits) of $X$ without sacrificing
information is still an open problem. We believe that fundamental weights and
their dual may add, but do not solve the quest. Finally, Section
\ref{sec_conclusion} briefly concludes.

\section{Eigenvector components as determinants}

\label{sec_determinantal_square_eigenvector_component}We assume that the
eigenvalue $\lambda_{k}$ is single with multiplicity one, so that rank$\left(
A-\lambda_{k}I\right)  =N-1$. This means that $\left(  A-\lambda_{k}I\right)
x_{k}=0$ contains only $N-1$ linearly independent equations to determine the
$N$ unknowns $\left(  x_{k}\right)  _{1},\left(  x_{k}\right)  _{2}%
,\ldots,\left(  x_{k}\right)  _{N}$. There are basically two
approaches\footnote{These two approaches are similar to computing the adjoint
matrix $Q\left(  \lambda\right)  =c_{A}\left(  \lambda\right)  \left(  \lambda
I-A\right)  ^{-1}$, whose columns are eigenvectors (see \cite[art. 148 on p.
220]{PVM_graphspectra}, \cite[Chapter IV]{GantmacherI}).} to determine the $N$
unknowns: (i) one of the $N$ equations/rows in $A-\lambda_{k}I$ can be
replaced by an additional equation as explored below and (ii) the set is
rewritten in $N-1$ unknowns in terms of one of them, say $\left(
x_{k}\right)  _{N}$, whose analysis is omitted, because the resulting
expressions for $\left(  x_{k}\right)  _{j}$ are less general as those in (i).

We replace an arbitrary equation or row in the set $\left(  A-\lambda
_{k}I\right)  x_{k}=0$ by a new linear equation $\,b^{T}x_{k}=\sum_{j=1}%
^{N}b_{j}\left(  x_{k}\right)  _{j}$, where $b$ is a real vector and the real
number $\beta_{k}=b^{T}x_{k}$ is non-zero. In most cases (except for regular
graphs where the all-one vector $u=\left(  1,1,\ldots,1\right)  $ is an
eigenvector), that additional equation is a normalization relation for the
eigenvector and the simplest linear one is $u^{T}x_{k}=w_{k}$, where
$w_{k}\neq0$ is a real number and called \emph{the fundamental weight}
\cite{PVM_double_orthogonality,PVM_Lower_bound_fundamental_weight_SITIS2014}
of $x_{k}$, further discussed in Section \ref{sec_fundamental_weights} while
formulas for $\beta_{k}$ are summarized in Appendix
\ref{sec_addition_Theorem_square_xk_determinant}. Another example is the
degree vector, $b=d$, where $d^{T}x_{k}=\lambda_{k}w_{k}$. The general
orthogonality equation $x_{k}^{T}x_{m}=\sum_{j=1}^{N}\left(  x_{k}\right)
_{j}\left(  x_{m}\right)  _{j}=\delta_{km}$ is another linear equation in the
unknown components of the vector $x_{k}$, given the components of the vector
$x_{m}$. However, since in this case $x_{k}^{T}x_{m}=0$, those linear
equations cannot be used!

\begin{theorem}
\label{theorem_eigenvector_component_determinants}Let $A$ and $A_{G\backslash
\left\{  j\right\}  }$ denote the adjacency matrix of the graph $G$ and of the
graph $G_{\backslash\left\{  j\right\}  }$ in which node $j$ and all its
incident links are removed from $G$, respectively. For any vector $b$ with
$\beta_{k}=b^{T}x_{k}\neq0$, the $j$-th component of eigenvector $x_{k}$ of
$A$ belonging to eigenvalue $\lambda_{k}$ can be written as
\begin{equation}
\left(  x_{k}\right)  _{j}=\frac{\beta_{k}\det\left(  A_{G\backslash\left\{
j\right\}  }-\lambda_{k}I\right)  }{\det\left(  A-\lambda_{k}I\right)
_{\operatorname{row}j=b}} \label{x_k_j_detA/j}%
\end{equation}
or%
\begin{equation}
\left(  x_{k}\right)  _{j}=-\frac{\det\left(  A-\lambda_{k}I\right)
_{\operatorname{row}j=b}}{\beta_{k}c_{A}^{\prime}\left(  \lambda_{k}\right)  }
\label{eigenvector_xk_component_j}%
\end{equation}
where $\det\left(  A-\lambda_{k}I\right)  _{\operatorname{row}j=b}$ is the
$N\times N$ matrix obtained from $\left(  A-\lambda_{k}I\right)  $ by
replacing row $j$ by the vector $b$. The square of the $j$-th component of
eigenvector $x_{k}$ of $A$ belonging to eigenvalue $\lambda_{k}$ with
multiplicity 1 equals
\begin{equation}
\left(  x_{k}\right)  _{j}^{2}=-\frac{1}{c_{A}^{\prime}\left(  \lambda
_{k}\right)  }\det\left(  A_{G\backslash\left\{  j\right\}  }-\lambda
_{k}I\right)  \label{squared_eigenvector_component_xk_j_charc_pol}%
\end{equation}
where $c_{A}\left(  \lambda\right)  =\det\left(  A-\lambda I\right)  $ is the
characteristic polynomial of $A$ and $c_{A}^{\prime}\left(  \lambda\right)
=\frac{dc_{A}\left(  \lambda\right)  }{d\lambda}$.
\end{theorem}

Although formulated in terms of the adjacency matrix of a graph, Theorem
\ref{theorem_eigenvector_component_determinants} holds for any symmetric matrix.

\textbf{Proof:} Without loss of generality, we first replace the $N$-th
equation in $\left(  A-\lambda_{k}I\right)  x_{k}=0$ by $b^{T}x_{k}=\beta_{k}$
and the resulting set of linear equations becomes%
\[
\left[
\begin{array}
[c]{c}%
\left(  A-\lambda_{k}I\right)  _{\backslash\operatorname{row}N}\\
b
\end{array}
\right]  x_{k}=\left[
\begin{array}
[c]{c}%
0_{\left(  N-1\right)  \times1}\\
\beta_{k}%
\end{array}
\right]
\]
where $\left(  A-\lambda_{k}I\right)  _{\backslash\operatorname{row}N}$ is the
$\left(  N-1\right)  \times N$ matrix obtained from $\left(  A-\lambda
_{k}I\right)  $ by removing row $N$. Cramer's solution \cite[p. 256]%
{PVM_graphspectra} yields%
\[
\left(  x_{k}\right)  _{j}=\frac{\left\vert
\begin{array}
[c]{c}%
\left(  A-\lambda_{k}I\right)  _{\backslash\operatorname{row}N}\\
b
\end{array}
\right\vert _{\operatorname{col}j=\left[
\begin{array}
[c]{c}%
0_{\left(  N-1\right)  \times1}\\
\beta_{k}%
\end{array}
\right]  }}{\left\vert
\begin{array}
[c]{c}%
\left(  A-\lambda_{k}I\right)  _{\backslash\operatorname{row}N}\\
b
\end{array}
\right\vert }=\frac{\left(  -1\right)  ^{N+j}\beta_{k}\det\left(
A-\lambda_{k}I\right)  _{\backslash\operatorname{row}N\backslash
\operatorname{col}j}}{\det\left(  A-\lambda_{k}I\right)  _{\operatorname{row}%
N=b}}%
\]
The $j$-th component of the $k$-th eigenvector $x_{k}$ can be written
as\footnote{Remark that the adjacency matrix $A_{G\backslash\operatorname{row}%
m\backslash\operatorname{col}i}$ represents a directed graph in which the
out-going links of node $m$ and the in-coming links to node $i$ are removed;
everywhere else, the in-coming and out-going links are the same
(bidirectional). Thus, $A_{G\backslash\operatorname{row}m\backslash
\operatorname{col}i}$ is not necessarily symmetric and it has $\left\vert
m-i\right\vert $ non-zero diagonal elements, $a_{k+1,k}$ for $m\leq k<i$.}%

\begin{equation}
\left(  x_{k}\right)  _{j}=\alpha_{m}\left(  k\right)  \left(  -1\right)
^{j}\det\left(  A-\lambda_{k}I\right)  _{\backslash\operatorname{row}%
m\backslash\operatorname{col}j} \label{eigenvector_component_xk_j}%
\end{equation}
where we have now deleted row $1\leq m\leq N$, instead of row $N$ as before,
and where the scaling factor is%
\begin{equation}
\alpha_{m}\left(  k\right)  =\frac{\left(  -1\right)  ^{m}\beta_{k}}%
{\det\left(  A-\lambda_{k}I\right)  _{\operatorname{row}m=b}}
\label{def_scaling_alfa_m}%
\end{equation}
Combining (\ref{eigenvector_component_xk_j}) with (\ref{def_scaling_alfa_m})
for $m=j$ leads to (\ref{x_k_j_detA/j}).

We now impose the orthogonality equation $x_{k}^{T}x_{k}=1$. It follows from
(\ref{eigenvector_component_xk_j}) that%
\[
\left(  x_{k}\right)  _{j}^{2}=\alpha_{m}^{2}\left(  k\right)  \left(
\det\left(  A-\lambda_{k}I\right)  _{\backslash\operatorname{row}%
m\backslash\operatorname{col}j}\right)  ^{2}%
\]
Invoking the identity%
\begin{equation}
\left(  \det\left(  A_{G\backslash\operatorname{row}m\backslash
\operatorname{col}j}-\lambda I\right)  \right)  ^{2}=\det\left(
A_{G\backslash\left\{  m\right\}  }-\lambda I\right)  \det\left(
A_{G\backslash\left\{  j\right\}  }-\lambda I\right)  -\det\left(
A_{G\backslash\left\{  m,j\right\}  }-\lambda I\right)  \det\left(
A_{G}-\lambda I\right)  \label{Jacobi_generalized_cofactor_theorem_k=2}%
\end{equation}
which can be deduced from Jacobi's famous theorem of 1833 (see e.g. \cite[p.
25]{Mirsky}), yields%
\begin{align}
\alpha_{m}^{-2}\left(  k\right)  \left(  x_{k}\right)  _{j}^{2}  &
=\lim_{\lambda\rightarrow\lambda_{k}}\det\left(  A_{G\backslash\left\{
m\right\}  }-\lambda I\right)  \det\left(  A_{G\backslash\left\{  j\right\}
}-\lambda I\right)  -\det\left(  A_{G\backslash\left\{  m,j\right\}  }-\lambda
I\right)  \det\left(  A_{G}-\lambda I\right) \nonumber\\
&  =\det\left(  A_{G\backslash\left\{  m\right\}  }-\lambda_{k}I\right)
\det\left(  A_{G\backslash\left\{  j\right\}  }-\lambda_{k}I\right)
\label{intermediate_hulp}%
\end{align}
The condition $x_{k}^{T}x_{k}=\sum_{n=1}^{N}\left(  x_{k}\right)  _{n}^{2}=1$
specifies $\alpha_{m}\left(  k\right)  $ as
\begin{equation}
\alpha_{m}^{-2}\left(  k\right)  =\det\left(  A_{G\backslash\left\{
m\right\}  }-\lambda_{k}I\right)  \sum_{n=1}^{N}\det\left(  A_{G\backslash
\left\{  n\right\}  }-\lambda_{k}I\right)  \label{alfa_m_squared_inverse}%
\end{equation}
We observe that there is a degree of freedom via the choice of $m$. Thus, for
$m=j$ in (\ref{eigenvector_component_xk_j}), we obtain from
(\ref{intermediate_hulp}) and (\ref{alfa_m_squared_inverse})%
\begin{equation}
\left(  x_{k}\right)  _{j}^{2}=\frac{\det\left(  A_{G\backslash\left\{
j\right\}  }-\lambda_{k}I\right)  }{\sum_{n=1}^{N}\det\left(  A_{G\backslash
\left\{  n\right\}  }-\lambda_{k}I\right)  }
\label{squared_eigenvector_component_xk_j}%
\end{equation}
that is independent of the choice of the vector $b$. Since \cite{Meyer_matrix}%
\begin{equation}
\sum_{n=1}^{N}\det\left(  A_{G\backslash\left\{  n\right\}  }-\lambda
I\right)  =-\frac{d}{d\lambda}\det\left(  A-\lambda I\right)  =-c_{A}^{\prime
}\left(  \lambda\right)  \label{derivative_characteristic_polynomial_A}%
\end{equation}
we arrive at (\ref{squared_eigenvector_component_xk_j_charc_pol}). Combining
(\ref{x_k_j_detA/j}) and (\ref{squared_eigenvector_component_xk_j_charc_pol})
yields\footnote{We remark that taking the derivative of both sides of
(\ref{eigenvector_xk_component_j}) with respect to $b_{m}$ results in
(\ref{x_k_j_detA/j}).} (\ref{eigenvector_xk_component_j}). \hfill
$\square\medskip$

\textbf{Another proof of (\ref{squared_eigenvector_component_xk_j_charc_pol}%
):} We start from the resolvent \cite[p. 244]{PVM_graphspectra} of a symmetric
matrix $A$%
\[
\left(  A-zI\right)  _{jj}^{-1}=\frac{\det\left(  A_{\backslash\left\{
j\right\}  }-zI\right)  }{\det\left(  A-zI\right)  }=\sum_{m=1}^{N}%
\frac{\left(  x_{m}\right)  _{j}^{2}}{\lambda_{m}-z}%
\]
from which, using $c_{A}\left(  \lambda\right)  =\det\left(  A-\lambda
I\right)  =\prod_{j=1}^{N}\left(  \lambda_{j}-\lambda\right)  $,
\begin{align*}
\det\left(  A_{\backslash\left\{  j\right\}  }-\lambda_{k}I\right)   &
=\sum_{m=1}^{N}\left(  x_{m}\right)  _{j}^{2}\lim_{z\rightarrow\lambda_{k}%
}\frac{\prod_{j=1}^{N}\left(  \lambda_{j}-z\right)  }{\lambda_{m}-z}\\
&  =\left(  x_{k}\right)  _{j}^{2}\prod_{j=1;j\neq k}^{N}\left(  \lambda
_{j}-\lambda_{k}\right)
\end{align*}
Invoking (\ref{afgeleide_characteristic_pol_in_lambda_m}) yields
(\ref{squared_eigenvector_component_xk_j_charc_pol}).\hfill$\square\medskip$

The second proof of (\ref{squared_eigenvector_component_xk_j_charc_pol}),
written as $x_{j}^{2}=\frac{P_{G-j}\left(  \lambda\right)  }{P_{G}^{\prime
}\left(  \lambda\right)  }$ where $P_{G}\left(  z\right)  =\det\left(
A-zI\right)  $, has appeared earlier in Cvetcovic \emph{et al.} \cite[Theorem
3.1]{Cvetcovic_LAA2007}, who referred to Hagos \cite{Hagos_LAA2002}, who in
turn mentioned that Mukherjee and Datta \cite{Mukherjee_1989} (using a
perturbation technique) and Li and Feng (only for the largest eigenvalue) have
preceded him. Hagos \cite{Hagos_LAA2002} mentioned rightly that
\textquotedblleft Eq. (\ref{squared_eigenvector_component_xk_j_charc_pol}) is
probably not as well known as it should be\textquotedblright, which may
justify why we have placed (\ref{squared_eigenvector_component_xk_j_charc_pol}%
) in the abstract as well. In addition, Hagos \cite{Hagos_LAA2002} has shown
that (in our notation)%
\[
\sum_{l=1}^{r_{k}}\left(  x_{l}\right)  _{j}^{2}=\frac{r_{k}}{c_{A}^{\prime
}\left(  \lambda_{k}\right)  }\det\left(  A_{G\backslash\left\{  j\right\}
}-\lambda_{k}I\right)
\]
where $\lambda_{k}$ is an eigenvalue with multiplicity $r_{k}$ and $x_{l}$ is
one of the $r_{k}$ orthogonal eigenvectors belonging to eigenvalue
$\lambda_{k}$.

\begin{corollary}
\label{corollary_eigenvalue_multiplicity_2}If $\lambda_{k}$ is an eigenvalue
of $A$ with multiplicity of two, then%
\begin{equation}
\left(  x_{k}\right)  _{j}^{2}=\frac{1}{c_{A}^{^{\prime\prime}}\left(
\lambda_{k}\right)  }\sum_{n=1;n\neq j}^{N-1}\det\left(  A_{G\backslash
\left\{  j,n\right\}  }-\lambda_{k}I\right)
\label{eigenvector_component_multiciplicity_2}%
\end{equation}

\end{corollary}

\textbf{Proof:} If $\lambda_{k}$ is an eigenvalue of $A$ with multiplicity of
two, then it holds that $c_{A}\left(  \lambda_{k}\right)  =c_{A}^{\prime
}\left(  \lambda_{k}\right)  =0$. Moreover,
(\ref{derivative_characteristic_polynomial_A}) and the fact that $\det\left(
A_{G\backslash\left\{  m\right\}  }-\lambda_{k}I\right)  $ must have the same
sign (see e.g. (\ref{ratio_kwadraad}) below), show that all $\det\left(
A_{G\backslash\left\{  m\right\}  }-\lambda_{k}I\right)  $ must vanish,
implying that $\lambda_{k}$ is then also an eigenvalue of all $A_{G\backslash
\left\{  m\right\}  }$, for each node $m$ removed from $G$. This observation
agrees with the Interlacing theorem \cite{PVM_graphspectra} that tells us that
all eigenvalues of $A_{G\backslash\left\{  m\right\}  }$ (for each $m$) are
lying in between the eigenvalues of $A$. If two eigenvalues of $A$ coincide
(e.g. $\lambda_{k}=\lambda_{k+1}$), the corresponding eigenvalue of each
$A_{G\backslash\left\{  m\right\}  }$, i.e. $\lambda_{k}\geq\lambda\left(
A_{G\backslash\left\{  m\right\}  }\right)  \geq\lambda_{k+1}$, is squeezed to
that same value $\lambda_{k}$. Applying de l'Hospital's rule,
\[
\left(  x_{k}\right)  _{j}^{2}=-\lim_{\lambda\rightarrow\lambda_{k}}\frac
{\det\left(  A_{G\backslash\left\{  j\right\}  }-\lambda I\right)  }%
{c_{A}^{\prime}\left(  \lambda\right)  }=-\lim_{\lambda\rightarrow\lambda_{k}%
}\frac{\frac{d}{d\lambda}\det\left(  A_{G\backslash\left\{  j\right\}
}-\lambda I\right)  }{c_{A}^{^{\prime\prime}}\left(  \lambda\right)  }%
\]
The derivative (\ref{derivative_characteristic_polynomial_A}) yields%
\[
\sum_{n=1;n\neq j}^{N-1}\det\left(  A_{G\backslash\left\{  j,n\right\}
}-\lambda I\right)  =-\frac{d}{d\lambda}\det\left(  A_{G\backslash\left\{
j\right\}  }-\lambda I\right)
\]
Combining these formulas, leads to
(\ref{eigenvector_component_multiciplicity_2}). \hfill$\square\medskip$

If $c_{A}^{\prime}\left(  \lambda_{k}\right)  =0$,
(\ref{eigenvector_component_multiciplicity_2}) reflects the effect of removing
all pair of nodes containing node $j$.

\begin{corollary}
\label{corollary_prod_eigenvector_components}The product of the $j$-th and
$m$-th component of eigenvector $x_{k}$ of $A$ belonging to eigenvalue
$\lambda_{k}$ with multiplicity 1 equals%
\begin{equation}
\left(  x_{k}\right)  _{j}\left(  x_{k}\right)  _{m}=\frac{\left(  -1\right)
^{j+m+1}}{c_{A}^{\prime}\left(  \lambda_{k}\right)  }\det\left(
A_{\backslash\operatorname{row}j\backslash\operatorname{col}m}-\lambda
_{k}I\right)  \label{product_eigenvector_(xk)_j_and_m}%
\end{equation}

\end{corollary}

\textbf{Proof:} We expand the determinant in (\ref{eigenvector_xk_component_j}%
) in the cofactors of row $j$ and obtain, with $\beta_{k}=\sum_{m=1}^{N}%
b_{m}\left(  x_{k}\right)  _{m}$,%
\[
\sum_{m=1}^{N}b_{m}\left(  x_{k}\right)  _{m}\left(  x_{k}\right)  _{j}%
=-\frac{\left(  -1\right)  ^{j}}{c_{A}^{\prime}\left(  \lambda_{k}\right)
}\sum_{m=1}^{N}\left(  -1\right)  ^{m}b_{m}\det\left(  A_{\backslash
\operatorname{row}j\backslash\operatorname{col}m}-\lambda_{k}I\right)
\]
Since this relation holds for any vector $b=\left(  b_{1},b_{2},\ldots
,b_{N}\right)  $, equating the corresponding coefficient $b_{m}$ at both sides
yields (\ref{product_eigenvector_(xk)_j_and_m}). \hfill$\square\medskip$

When $m=j$ in (\ref{product_eigenvector_(xk)_j_and_m}), we arrive again at
(\ref{squared_eigenvector_component_xk_j_charc_pol}). Hence,
(\ref{product_eigenvector_(xk)_j_and_m}) generalizes
(\ref{squared_eigenvector_component_xk_j_charc_pol}). The second orthogonality
relation (\ref{second_orthogonality_equations}) indicates that%
\[
\left(  -1\right)  ^{j+m+1}\sum_{k=1}^{N}\frac{\det\left(  A_{\backslash
\operatorname{row}j\backslash\operatorname{col}m}-\lambda_{k}I\right)  }%
{c_{A}^{\prime}\left(  \lambda_{k}\right)  }=\delta_{jm}%
\]

\subsection{Interpretations}

\label{sec_interpretations}Additional deductions from Theorem
\ref{theorem_eigenvector_component_determinants} are presented in Appendix
\ref{sec_addition_Theorem_square_xk_determinant}.

\medskip\medskip\refstepcounter{article}{\noindent\textbf{\thearticle. }%
}\ignorespaces\label{art_component_ratios}\emph{Component ratios}. We deduce
from (\ref{squared_eigenvector_component_xk_j_charc_pol}) that%
\begin{equation}
\frac{\left(  x_{k}\right)  _{j}^{2}}{\left(  x_{k}\right)  _{m}^{2}}%
=\frac{\det\left(  A_{G\backslash\left\{  j\right\}  }-\lambda_{k}I\right)
}{\det\left(  A_{G\backslash\left\{  m\right\}  }-\lambda_{k}I\right)  }%
=\frac{c_{A_{G\backslash\left\{  j\right\}  }}\left(  \lambda_{k}\right)
}{c_{A_{G\backslash\left\{  m\right\}  }}\left(  \lambda_{k}\right)  }
\label{ratio_kwadraad}%
\end{equation}
illustrating that $\det\left(  A_{G\backslash\left\{  j\right\}  }-\lambda
_{k}I\right)  $ and $\det\left(  A_{G\backslash\left\{  m\right\}  }%
-\lambda_{k}I\right)  $ have the same sign for any pair of nodes $\left(
j,m\right)  $ for a given frequency $\lambda_{k}$, but, by
(\ref{derivative_characteristic_polynomial_A}), opposite to the sign of
$c_{A}^{\prime}\left(  \lambda_{k}\right)  $.

It follows from (\ref{eigenvector_xk_component_j}) that%
\begin{equation}
\frac{\left(  x_{k}\right)  _{j}}{\left(  x_{k}\right)  _{m}}=\frac
{\det\left(  A-\lambda_{k}I\right)  _{\operatorname{row}j=b}}{\det\left(
A-\lambda_{k}I\right)  _{\operatorname{row}m=b}}
\label{ration_eigenvector_components}%
\end{equation}
The sign of $\left(  x_{k}\right)  _{j}$ with respect to $\left(
x_{k}\right)  _{m}$ is thus determined by a ratio of determinants that
seemingly depend on an arbitrary vector $b$ with non-zero $\beta_{k}%
=u^{T}x_{k}$, whose general graph interpretation is less transparent than
nodal removal as in $\det\left(  A_{G\backslash\left\{  j\right\}  }%
-\lambda_{k}I\right)  $, even if $b=u$. If $k=1$, then $\left(  x_{1}\right)
_{j}\geq0$, so that $\det\left(  A-\lambda_{1}I\right)  _{\operatorname{row}%
j=b}$ and $\det\left(  A-\lambda_{1}I\right)  _{\operatorname{row}m=b}$ have
the same sign. However, for $k>1$, it holds that $\min_{j}\left(
x_{k}\right)  _{j}\leq0\leq\max_{j}\left(  x_{k}\right)  _{j}$ and, hence
(\ref{x_k_j_detA/j}) shows that $\min_{1\leq j\leq N}\det\left(  A-\lambda
_{k}I\right)  _{\operatorname{row}j=b}$ has a sign opposite to $\max_{1\leq
j\leq N}\det\left(  A-\lambda_{k}I\right)  _{\operatorname{row}j=b}$. We
remark that the ratios (\ref{ratio_kwadraad}) and
(\ref{ration_eigenvector_components}) only hold at eigenfrequencies of $A$,
thus%
\begin{equation}
\frac{\det\left(  A_{G\backslash\left\{  j\right\}  }-\lambda I\right)  }%
{\det\left(  A_{G\backslash\left\{  m\right\}  }-\lambda I\right)  }=\left(
\frac{\det\left(  A-\lambda I\right)  _{\operatorname{row}j=b}}{\det\left(
A-\lambda I\right)  _{\operatorname{row}m=b}}\right)  ^{2}
\label{ratio_determinants_at_eigenvalue_lambda_k}%
\end{equation}
is correct only if $\lambda=\lambda_{k}$ for $1\leq k\leq N$.

\medskip\medskip\refstepcounter{article}{\noindent\textbf{\thearticle. }%
}\ignorespaces\label{art_zero_eigenvector_component} \emph{Zero eigenvector
component.}\textbf{ }If $\lambda_{k}$ is a single eigenvalue of $A$ (thus
$c_{A}^{\prime}\left(  \lambda_{k}\right)  \neq0)$ and if $\lambda_{k}$ is
also an eigenvalue of $A_{G\backslash\left\{  j\right\}  }$, then
(\ref{squared_eigenvector_component_xk_j_charc_pol}) shows that $\left(
x_{k}\right)  _{j}=0$. Not all other eigenvector components $\left(
x_{k}\right)  _{m}$ can be zero, because any eigenvector is different from the
zero vector. Hence, if $\lambda_{k}$ is not an eigenvalue of multiplicity at
least two, then $\lambda_{k}$ cannot be an eigenvalue of all $A_{G\backslash
\left\{  m\right\}  }$ (for $1\leq m\leq N$). The eigenvalue equation states
that%
\[
\lambda_{k}\left(  x_{k}\right)  _{j}=\sum_{l=1}^{N}a_{jl}\left(
x_{k}\right)  _{l}=\sum_{l\in\text{ }\mathcal{N}_{j}}\left(  x_{k}\right)
_{l}%
\]
where $\mathcal{N}_{j}$ represents the set of direct neighbors of node $j$. A
zero eigenvector component, $\left(  x_{k}\right)  _{j}=0$ at eigenvalue
$\lambda_{k}$, means that (a) the average of the eigenvector components of the
neighbors of node $j$ is zero and (b) that node $j$ does not affect the
eigenvector component of any of its neighbors. When $\left(  x_{k}\right)
_{j}=0$, the removal of node $j$ has no effect at frequency $\lambda
_{k}<\lambda_{1}$. Since $\left(  x_{1}\right)  _{j}>0$ in a connected graph
(by the Perron-Frobenius Theorem), the removal of a node $j$ has always an
effect at eigenfrequency $\lambda_{1}$. Based on this notion, we may define
\emph{the redundancy }$r_{j}$ $\in\left[  0,N-1\right]  $ of node $j$ as the
number of eigenfrequencies at which $\left(  x_{k}\right)  _{j}=\left(
x_{k}\right)  _{j}^{2}=0$.

\begin{figure}
[h]
\begin{center}
\includegraphics[
height=10.7547cm,
width=16.1298cm
]%
{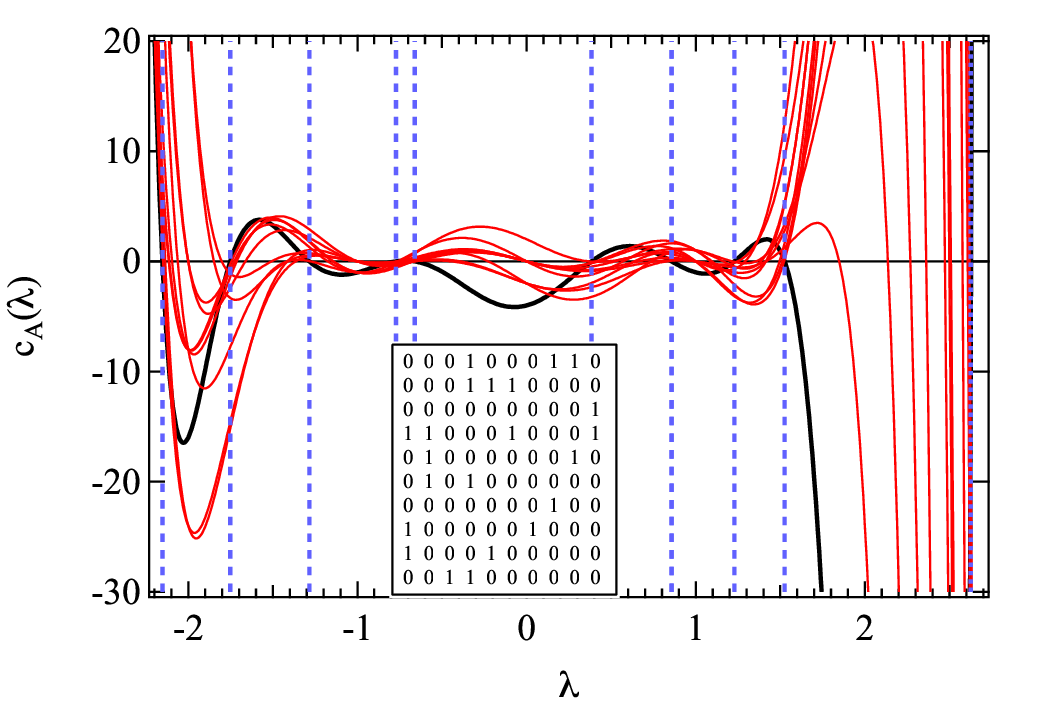}%
\caption{The characteristic polynomials $c_{A_{G\backslash\left\{  n\right\}
}}\left(  \lambda\right)  $ for $1\leq i\leq N$ in red and $c_{A_{G}}\left(
\lambda\right)  $ in black for an Erd\H{o}s-R\'{e}nyi graph $G_{0.2}\left(
10\right)  $, whose adjacency matrix is also shown. The blue vertical lines
denote the eigenvalues of $A$ (zeros of $c_{A}\left(  \lambda\right)  $).}%
\label{Fig_polynomialsERN10p02}%
\end{center}
\end{figure}

\medskip\medskip\refstepcounter{article}{\noindent\textbf{\thearticle. }%
}\ignorespaces\label{art_amplitude_interpretation} \emph{Amplitude.} The
magnitude of $\left(  x_{k}\right)  _{j}^{2}$ for node $j$ in
(\ref{squared_eigenvector_component_xk_j_charc_pol}) depends on the
characteristic polynomial $c_{A_{G\backslash\left\{  j\right\}  }}\left(
\lambda\right)  $ of $G\backslash\left\{  j\right\}  $ at the frequency
$\lambda=\lambda_{k}$. As illustrated in Fig.~\ref{Fig_polynomialsERN10p02},
the characteristic polynomials $c_{A}\left(  x\right)  $ and
$c_{A_{G\backslash\left\{  j\right\}  }}\left(  x\right)  $ oscillate around
zero in the interval $x\in\lbrack\lambda_{N},\lambda_{1}]$, that contains all
their real zeros. We coin the deviations in $c_{A_{G\backslash\left\{
j\right\}  }}\left(  x\right)  $ from zero at $\lambda_{k}$ the amplitude.
Just as in quantum mechanics (see e.g. \cite{Dirac,CohenTannoudji}), where the
wave function can be complex, while its modulus is interpreted as a
probability, we propose to use the eigenvector components $\left(
x_{k}\right)  _{j}$ in computations, but we suggest, based on
(\ref{squared_eigenvector_component_xk_j_charc_pol}), to interpret $\left(
x_{k}\right)  _{j}^{2}$ as centrality metrics. Hence, the importance or
centrality of node $j$ for property $\mathcal{P}_{k}$ at eigenfrequency
$\lambda_{k}$ is proportional to the amplitude of the characteristic
polynomial at $\lambda_{k}$ of the graph in which that node $j$ is removed.
Thus, the centrality $\left(  x_{k}\right)  _{j}^{2}$ measures a kind of
\textquotedblleft robustness\textquotedblright\ or \textquotedblleft
resilience\textquotedblright, in the sense of how important is the removal of
node $j$ from the graph $G$, determined by the amplitude at frequency
$\lambda_{k}$. In network robustness analyses, the removal of links or nodes
challenges the functioning of the network, measured via certain network
metrics
\cite{PVM_topologicalRobustness_evalutation,Manzano_Scoglio_NatScRe2014}. The
relative impact or effect of the removal of a high degree node at the largest
eigenfrequency $\lambda_{1}$ is larger than the removal of a low degree node
\cite{PVM_decreasingspectralradius_PRE2011}. However, at other
eigenfrequencies, the reverse must hold due to double orthogonality
(\ref{second_orthogonality_equations}), $\sum_{k=1}^{N}\left(  x_{k}\right)
_{j}^{2}=1$.

\begin{figure}
[h]
\begin{center}
\includegraphics[
height=10.7547cm,
width=16.1298cm
]%
{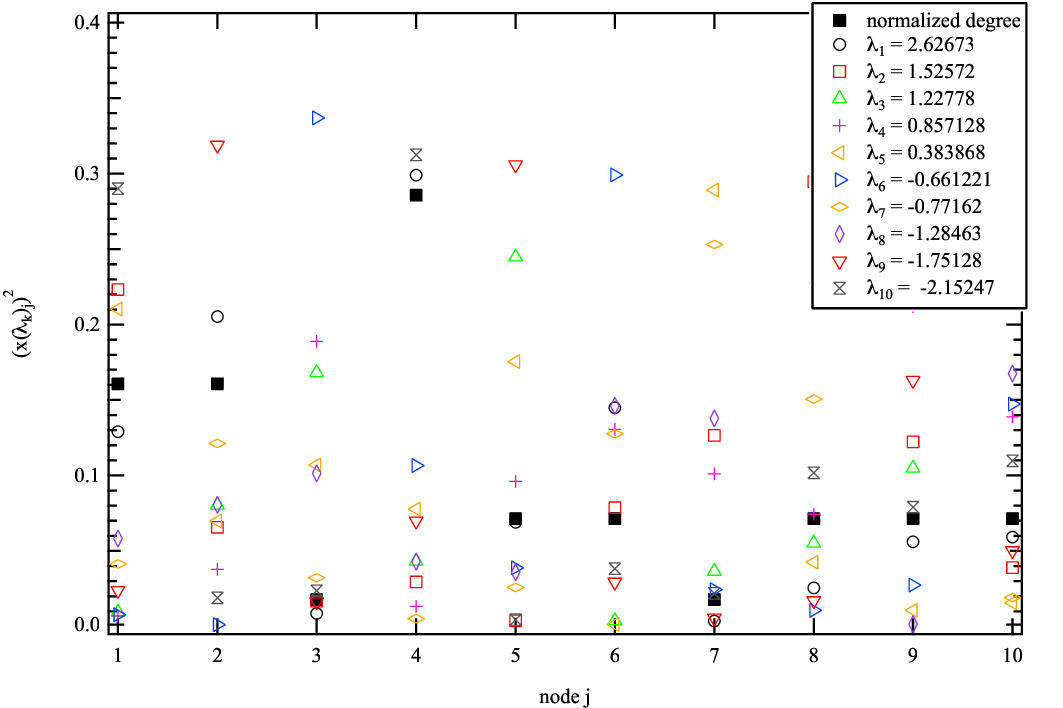}%
\caption{The square of the eigenvector components per node $j$ over all
eigenvalues $\lambda_{k}$ for the same graph as in
Fig.~\ref{Fig_polynomialsERN10p02}. The filled black squares represent the
normalized degree $d_{j}^{2}/d^{T}d$.}%
\label{Fig_squaredeigenvectorcomponentsern10p02}%
\end{center}
\end{figure}

\medskip\textbf{Example.} For a connected Erd\H{o}s-R\'{e}nyi graph with link
density $p=0.2$, $N=10$ nodes and the degree vector $d=\left(
3,3,1,4,2,2,1,2,2,2\right)  $, Fig.~\ref{Fig_polynomialsERN10p02} shows all 10
characteristic polynomials\footnote{The explicit expressions are%
\begin{align*}
c_{A}\left(  x\right)   &  =-\ 4+4x+27x^{2}-10x^{3}-52x^{4}+8x^{5}%
+38x^{6}-2x^{7}-11x^{8}+x\ ^{10}\\
c_{A_{G\backslash\left\{  1\right\}  }}\left(  x\right)   &  =-2-5x+6x^{2}%
+17x^{3}-6x^{4}-19x^{5}+2x^{6}+8x^{7}-x\ ^{9}\\
c_{A_{G\backslash\left\{  2\right\}  }}\left(  x\right)   &  =-4x+16x^{3}%
-19x^{5}+8x^{7}-x^{9}\\
c_{A_{G\backslash\left\{  3\right\}  }}\left(  x\right)   &  =-8x+4x^{2}%
+29x^{3}-6x^{4}-29x^{5}+2x^{6}+10x^{7}-x^{9}\\
c_{A_{G\backslash\left\{  4\right\}  }}\left(  x\right)   &  =-4x+14x^{3}%
-16x^{5}+7x^{7}-x^{9}\\
c_{A_{G\backslash\left\{  5\right\}  }}\left(  x\right)   &  =-2-5x+8x^{2}%
+20x^{3}-8x^{4}-23x^{5}+2x^{6}+9x^{7}-x\ ^{9}\\
c_{A_{G\backslash\left\{  6\right\}  }}\left(  x\right)   &  =2-7x-4x^{2}%
+25x^{3}+2x^{4}-25x^{5}+9x^{7}-x^{9}\\
c_{A_{G\backslash\left\{  7\right\}  }}\left(  x\right)   &  =-2-9x+6x^{2}%
+30x^{3}-6x^{4}-29x^{5}+2x^{6}+10x^{7}-x\ ^{9}\\
c_{A_{G\backslash\left\{  8\right\}  }}\left(  x\right)   &  =-4x+2x^{2}%
+18x^{3}-4x^{4}-22x^{5}+2x^{6}+9x^{7}-x^{9}\\
c_{A_{G\backslash\left\{  9\right\}  }}\left(  x\right)   &  =-4x+4x^{2}%
+20x^{3}-6x^{4}-23x^{5}+2x^{6}+9x^{7}-x^{9}\\
c_{A_{G\backslash\left\{  10\right\}  }}\left(  x\right)   &  =-4x+4x^{2}%
+19x^{3}-6x^{4}-23x^{5}+2x^{6}+9x^{7}-x^{9}%
\end{align*}
} $c_{A_{G\backslash\left\{  j\right\}  }}\left(  \lambda\right)  $ and
$c_{A}\left(  \lambda\right)  $, as well as its adjacency matrix $A$. At the
vertical lines, that indicate the positions of the eigenvalues of $A$, all
values $c_{A_{G\backslash\left\{  j\right\}  }}\left(  \lambda_{k}\right)  $
for $1\leq j\leq10$ have a same sign, in agreement with (\ref{ratio_kwadraad}%
). The amplitude $c_{A_{G\backslash\left\{  j\right\}  }}\left(  \lambda
_{k}\right)  $ is a relative measure for $\left(  x_{k}\right)  _{j}^{2}$ and
indicates the importance of node $j$ at frequency $\lambda_{k}$.
Fig.~\ref{Fig_squaredeigenvectorcomponentsern10p02} illustrates that the
topological degree vector $d$ correlates best with the square components of
the principal eigenvector $x_{1}$. At other eigenfrequencies, other nodes are
\textquotedblleft important\textquotedblright.
Fig.~\ref{Fig_squaredeigenvectorcomponentsern10p02} also shows that $\left(
x_{1}\right)  _{j}^{2}=\min_{1\leq k\leq\leq10}\left(  x_{k}\right)  _{j}^{2}$
for node $j=3$ and $j=7$, both having the minimum degree $d_{\min}=1$.

\medskip\refstepcounter{article}{\noindent\textbf{\thearticle. }%
}\ignorespaces\label{art_lack_independence_centrality_metrics}\emph{Concern
for the adjacency matrix }$A$: The zero eigenvalue in
(\ref{zero_eigenvalue_equation_Y}) of $C$ implies for any adjacency matrix $A$
that rank$\left(  C\right)  <N$ and that at least one row (or column) is a
linear combination of all the other rows (columns). Hence, the set of
centrality metrics $\left\{  \left(  \operatorname{row}C\right)  _{i}\right\}
_{1\leq i\leq N}=\left\{  \left(  x_{1}\right)  _{i}^{2},\left(  x_{2}\right)
_{i}^{2},\ldots,\left(  x_{N}\right)  _{i}^{2}\right\}  _{1\leq i\leq N}$ is
\emph{not} independent for the adjacency matrix, indicating that the set of
centrality metrics belonging to node $j$ can be written in terms of the
centrality metrics of all the others nodes in $G$.

\medskip\refstepcounter{article}{\noindent\textbf{\thearticle. }%
}\ignorespaces\label{art_link_addition_removal}\emph{Link addition/removal to
the graph }$G$\emph{.} Equation
(\ref{squared_eigenvector_component_xk_j_charc_pol}) indicates that the
addition (or removal) of a link to node $j$ does not change $\left(
x_{k}\right)  _{j}$, because $G_{\backslash\left\{  j\right\}  }$ means that,
besides the node $j$ itself, also all incident links to node $j$ are removed
from the graph. However, a link addition/removal may change the
eigenfrequencies $\left\{  \lambda_{k}\right\}  _{1\leq k\leq N}$. This
observation may suggest that, after the addition (or removal) of a link to
node $i$ and node $j$, the nodal eigenvector component $\left(  x_{k}\right)
_{i}$ and $\left(  x_{k}\right)  _{j}$ change the least. Simulations do not
seem to support this observation, which hints that the effect of link
addition/removal on the eigenfrequencies is dominant.

\medskip\refstepcounter{article}{\noindent\textbf{\thearticle. }%
}\ignorespaces\label{art_weighting_centralities}\emph{Weighting squared
eigenvector components.} Let $f\left(  x\right)  =x^{2}$ in
(\ref{eigenvalue_decompostion_function_f}), then%
\[
A^{2}=\sum_{k=1}^{N}\lambda_{k}^{2}x_{k}x_{k}^{T}=\sum_{k=1}^{N}\left(
\left\vert \lambda_{k}\right\vert x_{k}\right)  \left(  \left\vert \lambda
_{k}\right\vert x_{k}\right)  ^{T}%
\]
On the other hand, for the Laplacian $Q=\Delta-A$ whose eigenvalues are
non-negative, (\ref{eigenvalue_decompostion_function_f}) with $f\left(
x\right)  =x$ becomes%
\[
Q=\sum_{k=1}^{N}\mu_{k}z_{k}z_{k}^{T}=\sum_{k=1}^{N}\left(  \sqrt{\mu_{k}%
}z_{k}\right)  \left(  \sqrt{\mu_{k}}z_{k}\right)  ^{T}%
\]
These relations suggest to weight the \textquotedblleft
importance\textquotedblright\ of the eigenvectors of $A$ as $v_{k}=\left\vert
\lambda_{k}\right\vert x_{k}$, whereas those of $Q$ as $s_{k}=\sqrt{\mu_{k}%
}z_{k}$. Moreover, since $\left(  A^{2}\right)  _{jj}=Q_{jj}=d_{j}$ and
$\mu_{N}=0$, the two expression for the degree%
\[
d_{j}=\sum_{k=1}^{N}\lambda_{k}^{2}\left(  x_{k}^{2}\right)  _{j}=\sum
_{k=1}^{N-1}\mu_{k}\left(  z_{k}^{2}\right)  _{j}%
\]
show a weighting of the adjacency eigenvector centralities $\left(  x_{k}%
^{2}\right)  _{j}$ by $\lambda_{k}^{2}$, whereas the Laplacian eigenvector
centralities are only weighted proportional with the Laplacian eigenvalue
$\mu_{k}$. Thus, while the eigenvectors of different graph-related matrices
reflect different properties of the graph, although each of them satisfies the
first (\ref{first_orthogonality_equations}) and second
(\ref{second_orthogonality_equations}) orthogonality conditions, the example
illustrates that a generally acceptable scaling or weighting does not exist.
Clearly, the eigenvectors corresponding to the larger (in absolute value)
eigenvalues deserve more weight, as earlier was exploited in graph
reconstructability \cite{PVM_reconstructability_PRE} and only a few of the
larger ones may be sufficient as centrality metrics.

\section{Squared eigenvalue equation}

\label{sec_squared_eigenvalue_equation}

\begin{theorem}
\label{theorem_square_eigenvector_component_error_r}The square of the $i$-th
component of the eigenvector $x_{k}$ of the adjacency matrix $A$ of the graph
$G$ belonging to the eigenvalue $\lambda_{k}$ equals%
\begin{equation}
(x_{k})_{i}^{2}=\frac{1-r_{i}\left(  k\right)  }{\frac{\lambda_{k}^{2}%
(A)}{d_{i}}+1} \label{i_th_component_xk_square}%
\end{equation}
where $d_{i}$ is the degree of node $i$ and%
\begin{equation}
r_{i}\left(  k\right)  =\sum_{j=1;j\neq i}^{N}\left(  1-a_{ij}\right)
(x_{k})_{j}^{2}+\frac{1}{2d_{i}}\sum_{j=1}^{N}a_{ij}\sum_{l=1}^{N}%
a_{il}\left(  (x_{k})_{l}-(x_{k})_{j}\right)  ^{2} \label{def_r_i(k)}%
\end{equation}
obeys $0\leq r_{i}\left(  k\right)  \leq1$.
\end{theorem}

\textbf{Proof:} We start from the squared eigenvalue equation%
\[
\lambda_{k}^{2}(A)(x_{k})_{i}^{2}=\left(  \sum_{j=1}^{N}a_{ij}(x_{k}%
)_{j}\right)  ^{2}%
\]
to deduce an approximation for $(x_{k})_{i}^{2}$. Invoking the Cauchy identity
\cite[p. 257]{PVM_graphspectra} and $a_{ij}=a_{ij}^{2}$ yields%
\begin{align*}
\left(  \sum_{j=1}^{N}a_{ij}a_{ij}(x_{k})_{j}\right)  ^{2}  &  =\sum_{j=1}%
^{N}a_{ij}^{2}\sum_{j=1}^{N}\left(  a_{ij}(x_{k})_{j}\right)  ^{2}-\frac{1}%
{2}\sum_{j=1}^{N}\sum_{l=1}^{N}\left(  a_{ij}a_{il}(x_{k})_{l}-a_{il}%
a_{ij}(x_{k})_{j}\right)  ^{2}\\
&  =d_{i}\sum_{j=1}^{N}a_{ij}(x_{k})_{j}^{2}-\frac{1}{2}\sum_{j=1}^{N}%
a_{ij}\sum_{l=1}^{N}a_{il}\left(  (x_{k})_{l}-(x_{k})_{j}\right)  ^{2}%
\end{align*}
where the degree $d_{i}=\sum_{j=1}^{N}a_{ij}$. Further, using the first
orthogonality relations (\ref{first_orthogonality_equations}), $1=\sum
_{j=1}^{N}(x_{k})_{j}^{2}$, and%
\[
\sum_{j=1}^{N}a_{ij}(x_{k})_{j}^{2}=1-(x_{k})_{i}^{2}-\sum_{j=1;j\neq i}%
^{N}\left(  1-a_{ij}\right)  (x_{k})_{j}^{2}%
\]
we obtain%
\[
\frac{\lambda_{k}^{2}(A)}{d_{i}}(x_{k})_{i}^{2}=1-(x_{k})_{i}^{2}%
-\sum_{j=1;j\neq i}^{N}\left(  1-a_{ij}\right)  (x_{k})_{j}^{2}-\frac
{1}{2d_{i}}\sum_{j=1}^{N}a_{ij}\sum_{l=1}^{N}a_{il}\left(  (x_{k})_{l}%
-(x_{k})_{j}\right)  ^{2}%
\]
which we rewrite as (\ref{i_th_component_xk_square}). The definition
(\ref{def_r_i(k)}) shows that $r_{i}\left(  k\right)  \geq0$, whereas it
follows from (\ref{i_th_component_xk_square}) that $r_{i}\left(  k\right)
\leq1$.\hfill$\square\medskip$

Since $r_{i}\left(  k\right)  \geq0$, Theorem
\ref{theorem_square_eigenvector_component_error_r} directly leads to the upper
bound%
\begin{equation}
(x_{k})_{i}^{2}\leq\frac{1}{1+\frac{\lambda_{k}^{2}(A)}{d_{i}}}
\label{bound_improved_quadrate_eigenvector_component_xk_adjacency_matrix}%
\end{equation}
which appeared earlier for $k=1$ in \cite{Cioaba_Gregory_ELA2007} and \cite[p.
29]{Stevanovic2015}. Equality in
(\ref{bound_improved_quadrate_eigenvector_component_xk_adjacency_matrix}) only
holds if $r_{i}\left(  k\right)  =0$, which is equivalent to both%
\[
\sum_{j=1;j\neq i}^{N}\left(  1-a_{ij}\right)  (x_{k})_{j}^{2}=\sum
_{j\notin\mathcal{N}_{i}}^{N}(x_{k})_{j}^{2}=0
\]
where $\mathcal{N}_{i}$ is the set of all direct neighbors of node $i$, and
$\left(  x_{k}\right)  _{l}=\left(  x_{k}\right)  _{j}$ for all nodes
$l,j\in\mathcal{N}_{i}$. In conclusion, for any eigenfrequency $k$, equality
in (\ref{bound_improved_quadrate_eigenvector_component_xk_adjacency_matrix})
is only possible if $\left(  x_{k}\right)  _{j}=0$ for $j\notin\mathcal{N}%
_{i}$ and $\left(  x_{k}\right)  _{l}=\frac{\pm1}{\sqrt{d_{i}}}$ for
$l\in\mathcal{N}_{i}$. If $k=1$, equality can only happen in a disconnected
graph consisting of a regular graph on $d_{i}$ nodes (thus the complete graph
$K_{d_{i}}$) and $N-d_{i}$ disconnected nodes from node $i$.

\subsection{Bounds for eigenvector components}

We present a number of bounds for the minimum and maximum of eigenvector
components, either over frequencies $k$ or over nodes $j$.

We remark as in \cite[p. 31]{Stevanovic2015} that $\min_{1\leq j\leq N}%
(x_{k})_{j}^{2}$, deduced from
(\ref{bound_improved_quadrate_eigenvector_component_xk_adjacency_matrix}), can
be sharpened.

\begin{corollary}
\label{corollary_upper_bound_min_(x_k)_j}For any graph, it holds that%
\begin{equation}
\min_{1\leq j\leq N}(x_{k})_{j}^{2}\leq\frac{1-\frac{d_{\min}}{2}s_{k}}%
{\frac{\lambda_{k}^{2}(A)}{d_{\min}}+N-d_{\min}}
\label{minimum_component_Nikiforov}%
\end{equation}
where $s_{k}=\min_{l,j}\left(  (x_{k})_{l}-(x_{k})_{j}\right)  ^{2}$ is the
minimal square spacing between eigenvector components of $x_{k}$.
\end{corollary}

\textbf{Proof}: The definition (\ref{def_r_i(k)}) of $r_{i}\left(  k\right)  $
reveals that%
\[
\sum_{j=1;j\neq i}^{N}\left(  1-a_{ij}\right)  (x_{k})_{j}^{2}\geq\left(
N-1-d_{i}\right)  \min_{1\leq j\leq N}(x_{k})_{j}^{2}%
\]
and%
\[
\sum_{j=1}^{N}a_{ij}\sum_{l=1}^{N}a_{il}\left(  (x_{k})_{l}-(x_{k}%
)_{j}\right)  ^{2}\geq d_{i}^{2}\min_{l,j}\left(  (x_{k})_{l}-(x_{k}%
)_{j}\right)  ^{2}=d_{i}^{2}s_{k}%
\]
so that%
\[
r_{i}\left(  k\right)  \geq\left(  N-1-d_{i}\right)  \min_{1\leq j\leq
N}(x_{k})_{j}^{2}+\frac{d_{i}}{2}s_{k}%
\]
Hence, (\ref{i_th_component_xk_square}) can be bounded
\[
\min_{1\leq j\leq N}(x_{k})_{j}^{2}\leq(x_{k})_{i}^{2}=\frac{1-r_{i}\left(
k\right)  }{\frac{\lambda_{k}^{2}(A)}{d_{i}}+1}\leq\frac{1-\left(
N-1-d_{i}\right)  \min_{1\leq j\leq N}(x_{k})_{j}^{2}-\frac{d_{i}}{2}s_{k}%
}{\frac{\lambda_{k}^{2}(A)}{d_{i}}+1}%
\]
which holds for all nodes $i$, also for the node with minimum degree, leading
to (\ref{minimum_component_Nikiforov}).\hfill$\square\medskip$

Inequality (\ref{minimum_component_Nikiforov}) extends the result of Nikiforov
\cite{Nikiforov_odd_cycles_LAA2007} (where $k=1$ and the minimal square
spacing $s_{k}$ $=0$) to all eigenfrequencies $k$. The right-hand side of
(\ref{minimum_component_Nikiforov}) (with $s_{k}=0$) is minimized for $k=1$.
Since $\frac{\lambda_{k}^{2}(A)}{d_{i}}+1$ is maximal if $d_{i}=d_{\min}$ and
$k=1$, (\ref{i_th_component_xk_square}) shows that $\min_{1\leq k\leq N}%
(x_{k})_{i}^{2}$ is reached when $k=1$ at a minimum degree node if
$\max_{1\leq k\leq N}r_{i}\left(  k\right)  =r_{i}\left(  1\right)  $.
However, the minimum degree node $q$ does not always obey $\min_{1\leq k\leq
N}(x_{k})_{q}^{2}=\left(  x_{1}\right)  _{q}^{2}$.

Inspired by Cioab\u{a} and Gregory, we extend their Theorem 3.4 in
\cite{Cioaba_Gregory_ELA2007}:

\begin{theorem}
\label{theorem_bound_for_minimum_max_eigenvector_component}For any graph $G$,
it holds that%
\begin{equation}
\min_{1\leq j\leq N}(x_{k})_{j}\leq\frac{\lambda_{k}^{m}(A)w_{k}}{N_{m}}%
\leq\max_{1\leq j\leq N}(x_{k})_{j} \label{min_max_bound_(x_k)_j}%
\end{equation}
where $w_{k}=\sum_{j=1}^{N}\left(  x_{k}\right)  _{j}$ is fundamental weight
(\ref{def_fundamental_weight_wk}) and $N_{m}=u^{T}A^{m}u=\sum_{i=1}^{N}%
\sum_{j=1}^{N}\left(  A^{m}\right)  _{ij}$ is the total number of walks with
$m$ hops in the graph $G$. Furthermore, we have%
\begin{equation}
\frac{\left\vert \lambda_{k}^{m}(A)\right\vert }{\sqrt{N_{2m}}}\leq\max_{1\leq
j\leq N}(x_{k})_{j} \label{lower_bound_max_eigenvector_component}%
\end{equation}
The companion of (\ref{min_max_bound_(x_k)_j}) over frequencies $k$ is%
\begin{equation}
\min_{1\leq k\leq N}(x_{k})_{j}^{2}\leq\frac{\left(  A^{m}\right)  _{jj}%
}{W_{m}}\leq\max_{1\leq k\leq N}(x_{k})_{j}^{2} \label{min_max_bound_(x_k)_k}%
\end{equation}
where $W_{m}=\sum_{l=1}^{N}\left(  A^{m}\right)  _{ll}=\sum_{k=1}^{N}%
\lambda_{k}^{m}$ is the total number of closed walks \cite{PVM_graphspectra}
with $m$ hops/links.
\end{theorem}

\textbf{Proof:} Consider the eigenvalue equation%
\[
\lambda_{k}(A^{m})(x_{k})_{i}=\sum_{j=1}^{N}\left(  A^{m}\right)  _{ij}%
(x_{k})_{j}%
\]
First, we bound the sum%
\begin{equation}
\min_{1\leq j\leq N}(x_{k})_{j}\sum_{j=1}^{N}\left(  A^{m}\right)  _{ij}%
\leq\sum_{j=1}^{N}\left(  A^{m}\right)  _{ij}(x_{k})_{j}\leq\max_{1\leq j\leq
N}(x_{k})_{j}\sum_{j=1}^{N}\left(  A^{m}\right)  _{ij}
\label{min_max_bound_eigenvalue_equation_A^m}%
\end{equation}
and after introducing the above eigenvalue equation, we sum over all $i$ using
the definition (\ref{def_fundamental_weight_wk}) of the fundamental weight
$w_{k}$,%
\[
\min_{1\leq j\leq N}(x_{k})_{j}\sum_{i=1}^{N}\sum_{j=1}^{N}\left(
A^{m}\right)  _{ij}\leq\lambda_{k}(A^{m})w_{k}\leq\max_{1\leq j\leq N}%
(x_{k})_{j}\sum_{i=1}^{N}\sum_{j=1}^{N}\left(  A^{m}\right)  _{ij}%
\]
from which we find (\ref{min_max_bound_(x_k)_j}) with $\lambda_{k}%
(A^{m})=\lambda_{k}^{m}(A)$. Next, we square the inequality
(\ref{min_max_bound_eigenvalue_equation_A^m})%
\[
\lambda_{k}^{2}(A^{m})(x_{k})_{i}^{2}\leq\left(  \max_{1\leq j\leq N}%
(x_{k})_{j}\right)  ^{2}\sum_{j=1}^{N}\sum_{l=1}^{N}\left(  A^{m}\right)
_{ij}\left(  A^{m}\right)  _{il}%
\]
and then we sum over all $i$, using $\sum_{i=1}^{N}(x_{k})_{i}^{2}=1$,
\[
\lambda_{k}^{2}(A^{m})\leq\left(  \max_{1\leq j\leq N}(x_{k})_{j}\right)
^{2}\sum_{j=1}^{N}\sum_{l=1}^{N}\left(  \sum_{i=1}^{N}\left(  A^{m}\right)
_{ji}\left(  A^{m}\right)  _{il}\right)  =\left(  \max_{1\leq j\leq N}%
(x_{k})_{j}\right)  ^{2}\sum_{j=1}^{N}\sum_{l=1}^{N}\left(  A^{2m}\right)
_{jl}%
\]
which is equivalent to (\ref{lower_bound_max_eigenvector_component}).

For any non-negative function $f$, it follows directly from the general
formula (\ref{eigenvalue_decompostion_function_f_node_jj}) that%
\[
\min_{1\leq k\leq N}(x_{k})_{j}^{2}\leq\frac{\left(  f\left(  A\right)
\right)  _{jj}}{\sum_{k=1}^{N}f\left(  \lambda_{k}\right)  }\leq\max_{1\leq
k\leq N}(x_{k})_{j}^{2}%
\]
where $\sum_{k=1}^{N}f\left(  \lambda_{k}\right)  =\sum_{l=1}^{N}\left(
f\left(  A\right)  \right)  _{ll}$ (obtained by summing
(\ref{eigenvalue_decompostion_function_f_node_jj}) over all $j$ and invoking
(\ref{first_orthogonality_equations})). When choosing $f\left(  x\right)
=x^{m}$, we obtain (\ref{min_max_bound_(x_k)_k}).\hfill$\square\medskip$

The bound (\ref{min_max_bound_(x_k)_k}) illustrates that \textquotedblleft%
\ importance\textquotedblright\ of node $j$ over all eigenfrequencies $k$ is
dictated by the percentage of closed walks $\frac{\left(  A^{m}\right)  _{jj}%
}{W_{m}}$ of any length $m$ from and to that node $j$, which agrees with the
intuitive notion of importance in a network. For $m=2$,
(\ref{min_max_bound_(x_k)_k}) reduces with $W_{2}=2L=Nd_{av}$, where
$d_{av}=\frac{2L}{N}$ is the average degree in the graph $G$, to%
\[
\min_{1\leq k\leq N}(x_{k})_{j}^{2}\leq\frac{1}{N}\frac{d_{j}}{d_{av}}\leq
\max_{1\leq k\leq N}(x_{k})_{j}^{2}%
\]
while the case $m=0$ yields%
\[
\min_{1\leq k\leq N}(x_{k})_{j}^{2}\leq\frac{1}{N}\leq\max_{1\leq k\leq
N}(x_{k})_{j}^{2}%
\]
which illustrates that equality in both sides in (\ref{min_max_bound_(x_k)_k})
for irregular graphs is not possible. It follows from
(\ref{eigenvalue_decompostion_function_f_node_jj}) that $\left(  A^{m}\right)
_{jj}=\sum_{k=1}^{N}\lambda_{k}^{m}\left(  x_{k}\right)  _{j}^{2}$ so that,
for large $m$, $\left(  A^{m}\right)  _{jj}\sim\lambda_{1}^{m}\left(
x_{1}\right)  _{j}^{2}$ and $\sum_{k=1}^{N}\lambda_{k}^{m}\sim\lambda_{1}^{m}%
$, if $\lambda_{1}>\max\left(  \lambda_{2},\left\vert \lambda_{N}\right\vert
\right)  $, while $\left(  A^{m}\right)  _{jj}\sim2\lambda_{1}^{m}\left(
x_{1}\right)  _{j}^{2}$ and $\sum_{k=1}^{N}\lambda_{k}^{m}\sim2\lambda_{1}%
^{m}$ for complete bipartite graphs. Hence,%
\[
\lim_{m\rightarrow\infty}\frac{\left(  A^{m}\right)  _{jj}}{\sum_{k=1}%
^{N}\lambda_{k}^{m}}=\left(  x_{1}\right)  _{j}^{2}%
\]
and, for $m\rightarrow\infty$, the inequality (\ref{min_max_bound_(x_k)_k})
becomes $\min_{1\leq k\leq N}(x_{k})_{j}^{2}\leq\left(  x_{1}\right)  _{j}%
^{2}\leq\max_{1\leq k\leq N}(x_{k})_{j}^{2}$. Thus, the principal eigenvector
component can, in absolute value, be the smallest as well as the largest for a
node $j$ (see e.g. Fig.~\ref{Fig_squaredeigenvectorcomponentsern10p02}).

Combining (\ref{min_max_bound_(x_k)_j}) and
(\ref{lower_bound_max_eigenvector_component}) leads to%
\[
\max\left(  \frac{\lambda_{k}^{m}(A)w_{k}}{N_{m}},\frac{\left\vert \lambda
_{k}^{m}(A)\right\vert }{\sqrt{N_{2m}}}\right)  \leq\max_{1\leq j\leq N}%
(x_{k})_{j}%
\]
If $\lambda_{k}(A^{m})w_{k}>0$, then the inequality $N_{m}^{2}\leq NN_{2m}$
(see e.g. \cite[p. 34]{PVM_graphspectra}) does not allow us to deduce the
largest of the two lower bounds.

We now present another lower bound over all eigenfrequencies $k$.

\begin{corollary}
\label{corollary_lower_bound_maxeigenvector_component}The correction factor
$1-r_{i}\left(  k\right)  $, defined in (\ref{def_r_i(k)}), obeys%
\begin{equation}
\sum_{k=1}^{N}\left(  1-r_{i}\left(  k\right)  \right)  =2
\label{sum_(1-r_i(k))_over_k=2}%
\end{equation}
Moreover, the maximum eigenvector centrality $\left(  x_{k}\right)  _{i}^{2}$
of node $i$ is never smaller than%
\begin{equation}
\frac{4}{N\left(  3+\frac{\left(  A^{4}\right)  _{ii}}{d_{i}^{2}}\right)
}\leq\max_{1\leq k\leq N}\left(  x_{k}\right)  _{i}^{2}
\label{lower_bound_max_square_eigenvector_component}%
\end{equation}

\end{corollary}

\textbf{Proof:} Combining (\ref{degree_node_j_square_eigenvector_components})
and (\ref{i_th_component_xk_square}) directly yields\footnote{Directly summing
the definition (\ref{def_r_i(k)}) gives
\begin{align*}
\sum_{k=1}^{N}r_{i}\left(  k\right)   &  =\sum_{j=1;j\neq i}^{N}\left(
1-a_{ij}\right)  \sum_{k=1}^{N}(x_{k})_{j}^{2}+\frac{1}{2d_{i}}\sum_{j=1}%
^{N}a_{ij}\sum_{l=1}^{N}a_{il}\sum_{k=1}^{N}\left(  (x_{k})_{l}-(x_{k}%
)_{j}\right)  ^{2}\\
&  =\sum_{j=1;j\neq i}^{N}\left(  1-a_{ij}\right)  +\frac{1}{d_{i}}\sum
_{j=1}^{N}a_{ij}\sum_{l=1}^{N}a_{il}\left(  1-\delta_{lj}\right)
\end{align*}
where the second orthogonality relation (\ref{second_orthogonality_equations})
has been invoked. Further, with $\sum_{j=1;j\neq i}^{N}\left(  1-a_{ij}%
\right)  =\left(  N-1-d_{i}\right)  $ and
\begin{align*}
\frac{1}{d_{i}}\sum_{j=1}^{N}a_{ij}\sum_{l=1}^{N}a_{il}\left(  1-\delta
_{lj}\right)   &  =\frac{1}{d_{i}}\sum_{j=1}^{N}a_{ij}\left\{  \sum_{l=1}%
^{N}a_{il}-a_{ij}\right\}  =\frac{1}{d_{i}}\sum_{j=1}^{N}a_{ij}d_{i}-\frac
{1}{d_{i}}\sum_{j=1}^{N}a_{ij}^{2}\\
&  =d_{i}-1
\end{align*}
we arrive at (\ref{sum_(1-r_i(k))_over_k=2}).} (\ref{sum_(1-r_i(k))_over_k=2}%
). Via this method, thus using (\ref{i_th_component_xk_square}) and
(\ref{sum_(1-r_i(k))_over_k=2}), the variance of the numbers $\left\{
1-r_{i}\left(  1\right)  ,\ldots,\left(  1-r_{i}\left(  N\right)  \right)
\right\}  $ equals%
\begin{align*}
\text{Var}\left[  \left(  1-r_{i}\left(  k\right)  \right)  \right]   &
=\frac{1}{N}\sum_{k=1}^{N}\left(  1-r_{i}\left(  k\right)  \right)
^{2}-\left(  \frac{1}{N}\sum_{k=1}^{N}\left(  1-r_{i}\left(  k\right)
\right)  \right)  ^{2}\\
&  =\frac{1}{N}\sum_{k=1}^{N}\left(  x_{k}\right)  _{i}^{4}\left(
\frac{\lambda_{k}^{2}(A)}{d_{i}}+1\right)  ^{2}-\frac{4}{N^{2}}%
\end{align*}
The first term equals%
\[
\sum_{k=1}^{N}\left(  1-r_{i}\left(  k\right)  \right)  ^{2}=\sum_{k=1}%
^{N}\left(  x_{k}\right)  _{i}^{4}+\frac{2}{d_{i}}\sum_{k=1}^{N}\left(
x_{k}\right)  _{i}^{4}\lambda_{k}^{2}(A)+\frac{1}{d_{i}^{2}}\sum_{k=1}%
^{N}\left(  x_{k}\right)  _{i}^{4}\lambda_{k}^{4}(A)
\]
Further, with (\ref{degree_node_j_square_eigenvector_components}),%
\[
\sum_{k=1}^{N}\left(  x_{k}\right)  _{i}^{4}\lambda_{k}^{2}(A)\leq\max_{1\leq
k\leq N}\left(  x_{k}\right)  _{i}^{2}\sum_{k=1}^{N}\left(  x_{k}\right)
_{i}^{2}\lambda_{k}^{2}(A)=\max_{1\leq k\leq N}\left(  x_{k}\right)  _{i}%
^{2}d_{i}%
\]
and, similarly,%
\[
\sum_{k=1}^{N}\left(  x_{k}\right)  _{i}^{4}\lambda_{k}^{4}(A)\leq\max_{1\leq
k\leq N}\left(  x_{k}\right)  _{i}^{2}\sum_{k=1}^{N}\left(  x_{k}\right)
_{i}^{2}\lambda_{k}^{4}(A)=\max_{1\leq k\leq N}\left(  x_{k}\right)  _{i}%
^{2}\left(  A^{4}\right)  _{ii}%
\]
where $\left(  A^{4}\right)  _{ii}$ is the number of closed walks with 4 hops
starting and ending at node $i$, results in an upper bound for the variance%
\[
\text{Var}\left[  \left(  1-r_{i}\left(  k\right)  \right)  \right]  \leq
\frac{1}{N}\left(  \max_{1\leq k\leq N}\left(  x_{k}\right)  _{i}^{2}\left(
3+\frac{\left(  A^{4}\right)  _{ii}}{d_{i}^{2}}\right)  -\frac{4}{N}\right)
\]
Since the variance is non-negative, we find the lower bound
(\ref{lower_bound_max_square_eigenvector_component}).\hfill$\square\medskip$

Equation (\ref{sum_(1-r_i(k))_over_k=2}) indicates that the average over the
frequencies $k$ is $E_{k}\left[  1-r_{i}\left(  k\right)  \right]  =\frac
{1}{N}\sum_{k=1}^{N}\left(  1-r_{i}\left(  k\right)  \right)  =\frac{2}{N}$ so
that, approximately, $(x_{k})_{i}^{2}\approx\frac{1}{N}\frac{2d_{i}}%
{\lambda_{k}^{2}(A)+d_{i}}$.

\begin{theorem}
\label{theorem_bound_for_minimum_max_eigenvector_component copy(1)}For any
graph $G$, it holds that%
\begin{equation}
\min_{1\leq k\leq N}(x_{k})_{i}^{2}\leq\frac{1}{N}\frac{\min\left(
1+\frac{d_{av}}{d_{i}},2\right)  }{\min_{1\leq k\leq N}\left(  \frac
{\lambda_{k}^{2}(A)}{d_{i}}+1\right)  }
\label{upper_bound_minsquare_(x_k)_i_over_k}%
\end{equation}
and%
\begin{equation}
\min_{1\leq i\leq N}(x_{k})_{i}^{2}\leq\frac{1}{N}\frac{1+\lambda_{k}%
^{2}(A)E\left[  \frac{1}{D}\right]  }{1+\frac{\lambda_{k}^{2}(A)}{d_{\max}}}
\label{upper_bound_minsquare_(x_k)_i_over_i}%
\end{equation}
where the harmonic mean of the degree\footnote{As in
\cite{PVM_PAComplexNetsCUP}, the degree of a randomly chosen node in the graph
is denoted by the random variable $D$.} is $E\left[  \frac{1}{D}\right]
=\frac{1}{N}\sum_{i=1}^{N}\frac{1}{d_{i}}$.
\end{theorem}

\textbf{Proof:} Summing (\ref{i_th_component_xk_square}) over all $k$ and
invoking the second orthogonality relation
(\ref{second_orthogonality_equations}) yields%
\begin{equation}
1=\sum_{k=1}^{N}\frac{1-r_{i}\left(  k\right)  }{\frac{\lambda_{k}^{2}%
(A)}{d_{i}}+1}\text{ for nodes }1\leq i\leq N \label{unity_sum_r_i(k)}%
\end{equation}
while, similarly, the sum over all $i$ gives%
\[
1=\sum_{i=1}^{N}\frac{1-r_{i}\left(  k\right)  }{\frac{\lambda_{k}^{2}%
(A)}{d_{i}}+1}\text{ for frequency indices }1\leq k\leq N
\]
from which we obtain%
\[
\min_{1\leq k\leq N}\left(  1-r_{i}\left(  k\right)  \right)  \leq\left(
\sum_{k=1}^{N}\frac{1}{\frac{\lambda_{k}^{2}(A)}{d_{i}}+1}\right)  ^{-1}%
\]
and%
\[
\min_{1\leq i\leq N}\left(  1-r_{i}\left(  k\right)  \right)  \leq\left(
\sum_{i=1}^{N}\frac{1}{\frac{\lambda_{k}^{2}(A)}{d_{i}}+1}\right)  ^{-1}%
\]
Invoking the harmonic, geometric and arithmetic mean inequality (for positive,
real $a_{k}$)
\begin{equation}
\frac{n}{\sum_{k=1}^{n}\frac{1}{a_{k}}}\leq\sqrt[n]{\prod_{k=1}^{n}a_{k}}%
\leq\frac{1}{n}\sum_{k=1}^{n}a_{k}
\label{harmonic_geometric_arithmetic_mean_inequality}%
\end{equation}
shows, using $\sum_{k=1}^{N}\left(  \frac{\lambda_{k}^{2}(A)}{d_{i}}+1\right)
=N+\frac{2L}{d_{i}}$, that%
\[
\left(  \sum_{k=1}^{N}\frac{1}{\frac{\lambda_{k}^{2}(A)}{d_{i}}+1}\right)
^{-1}\leq\frac{1}{N}\left(  1+\frac{d_{av}}{d_{i}}\right)
\]
so that%
\[
\min_{1\leq k\leq N}\left(  1-r_{i}\left(  k\right)  \right)  \leq\left(
\sum_{k=1}^{N}\frac{1}{\frac{\lambda_{k}^{2}(A)}{d_{i}}+1}\right)  ^{-1}%
\leq\frac{1}{N}\left(  1+\frac{d_{av}}{d_{i}}\right)
\]
which is sharper than $\min_{1\leq k\leq N}\left(  1-r_{i}\left(  k\right)
\right)  \leq\frac{2}{N}$ (deduced from (\ref{sum_(1-r_i(k))_over_k=2})) when
$d_{av}\leq d_{i}$. Hence,%
\[
\min_{1\leq k\leq N}\left(  1-r_{i}\left(  k\right)  \right)  \leq\frac{1}%
{N}\min\left(  1+\frac{d_{av}}{d_{i}},2\right)
\]
so that, with $\min_{1\leq k\leq N}\left\{  (x_{k})_{i}^{2}\left(
\frac{\lambda_{k}^{2}(A)}{d_{i}}+1\right)  \right\}  \geq\min_{1\leq k\leq
N}(x_{k})_{i}^{2}\min_{1\leq k\leq N}\left(  \frac{\lambda_{k}^{2}(A)}{d_{i}%
}+1\right)  $, we obtain (\ref{upper_bound_minsquare_(x_k)_i_over_k}).
Similarly (for the node index), using $\sum_{i=1}^{N}\left(  \frac{\lambda
_{k}^{2}(A)}{d_{i}}+1\right)  =N+\lambda_{k}^{2}(A)\sum_{i=1}^{N}\frac
{1}{d_{i}}$, we have%
\[
\min_{1\leq i\leq N}\left(  1-r_{i}\left(  k\right)  \right)  \leq\frac{1}%
{N}\left(  1+\frac{\lambda_{k}^{2}(A)}{N}\sum_{i=1}^{N}\frac{1}{d_{i}}\right)
\]
and (\ref{harmonic_geometric_arithmetic_mean_inequality}) leads to $\frac
{1}{d_{av}}\leq E\left[  \frac{1}{D}\right]  =\frac{1}{N}\sum_{i=1}^{N}%
\frac{1}{d_{i}}\leq\frac{1}{d_{\min}}$. Invoking
\[
\min_{1\leq i\leq N}\left(  1-r_{i}\left(  k\right)  \right)  =\min_{1\leq
i\leq N}\left\{  (x_{k})_{i}^{2}\left(  \frac{\lambda_{k}^{2}(A)}{d_{i}%
}+1\right)  \right\}  \geq\min_{1\leq i\leq N}(x_{k})_{i}^{2}\left(
\frac{\lambda_{k}^{2}(A)}{d_{\max}}+1\right)
\]
finally yields (\ref{upper_bound_minsquare_(x_k)_i_over_i}).\hfill
$\square\medskip$

We observe that (\ref{upper_bound_minsquare_(x_k)_i_over_i}) is better for
small $\lambda_{k}^{2}(A)$ than (\ref{minimum_component_Nikiforov}) ignoring
$s_{k}$, while the opposite holds for large $\lambda_{k}^{2}(A)$.

Further, we bound (\ref{unity_sum_r_i(k)}), using
(\ref{sum_(1-r_i(k))_over_k=2}),%
\[
\frac{2}{\max_{1\leq k\leq N}\left(  \frac{\lambda_{k}^{2}(A)}{d_{i}%
}+1\right)  }\leq\sum_{k=1}^{N}\frac{1-r_{i}\left(  k\right)  }{\frac
{\lambda_{k}^{2}(A)}{d_{i}}+1}\leq\frac{2}{\min_{1\leq k\leq N}\left(
\frac{\lambda_{k}^{2}(A)}{d_{i}}+1\right)  }%
\]
and find%
\[
\min_{1\leq k\leq N}\left(  \lambda_{k}^{2}(A)\right)  \leq d_{i}\leq
\max_{1\leq k\leq N}\left(  \lambda_{k}^{2}(A)\right)  =\lambda_{1}^{2}(A)
\]
Since this inequality holds for each node $i$, we retrieve the classical bound
$\lambda_{1}(A)\geq\sqrt{d_{\max}}$ (equality holds for the star), but also
\[
\min_{1\leq k\leq N}\left(  \lambda_{k}^{2}(A)\right)  \leq d_{\min}%
\]
which is reminiscent to the inequality $\mu_{N-1}\leq d_{\min}$ for the
algebraic connectivity\footnote{The algebraic connectivity
\cite{Fiedler1973,PVM_graphspectra} is the second smallest eigenvalue
$\mu_{N-1}$ of the Laplacian $Q=\Delta-A$. Both the Laplacian $Q$ and $A^{2}$
have the same diagonal elements $Q_{jj}=\left(  A^{2}\right)  _{jj}=d_{j}$.}
(excluding the complete graph) and which we sharpen:

\begin{theorem}
\label{theorem_bound_minimum_eigenvalue_A^2}In any connected graph, it holds
that%
\begin{equation}
\min_{1\leq k\leq N}\left(  \lambda_{k}^{2}(A)\right)  <d_{\min}
\label{bound_min_eigenvalue_A^2}%
\end{equation}

\end{theorem}

\textbf{Proof:} Let us denote the ordering in the eigenvalues as
$\lambda_{\left(  1\right)  }^{2}(A)\geq$ $\lambda_{\left(  2\right)  }%
^{2}(A)\geq\ldots\geq\lambda_{\left(  N\right)  }^{2}(A)$, where
$\lambda_{\left(  1\right)  }^{2}(A)=\lambda_{1}^{2}(A)$ and $\lambda_{\left(
N\right)  }^{2}(A)=\min_{1\leq k\leq N}\left(  \lambda_{k}^{2}(A)\right)  $
and we write the index $k^{\ast}$ being associated with $\lambda_{\left(
k\right)  }^{2}(A)$, the $k$-th largest eigenvalue of $A^{2}$. After applying
Abel summation to (\ref{unity_sum_r_i(k)}), we obtain%
\[
\sum_{k=1}^{N}\frac{1-r_{i}\left(  k\right)  }{\frac{\lambda_{k}^{2}(A)}%
{d_{i}}+1}=\sum_{k=1}^{N-1}\left\{  \sum_{l=1}^{k}\left(  1-r_{i}\left(
l^{\ast}\right)  \right)  \right\}  \left(  \frac{1}{\frac{\lambda_{\left(
k\right)  }^{2}(A)}{d_{i}}+1}-\frac{1}{\frac{\lambda_{\left(  k+1\right)
}^{2}(A)}{d_{i}}+1}\right)  +\frac{1}{\frac{\lambda_{\left(  N\right)  }%
^{2}(A)}{d_{i}}+1}\sum_{l=1}^{N}\left(  1-r_{i}\left(  l^{\ast}\right)
\right)
\]
Using (\ref{sum_(1-r_i(k))_over_k=2}) and (\ref{unity_sum_r_i(k)}) yields%
\[
1=-S_{i}+\frac{2}{\frac{\lambda_{\left(  N\right)  }^{2}(A)}{d_{i}}+1}%
\]
where
\begin{equation}
S_{i}=\sum_{k=1}^{N-1}\left\{  \sum_{l=1}^{k}\left(  1-r_{i}\left(  l^{\ast
}\right)  \right)  \right\}  \left(  \frac{1}{\frac{\lambda_{\left(
k+1\right)  }^{2}(A)}{d_{i}}+1}-\frac{1}{\frac{\lambda_{\left(  k\right)
}^{2}(A)}{d_{i}}+1}\right)  \label{def_Si}%
\end{equation}
which is non-negative (because each term in the $k$-sum is), $S_{i}\geq0$.
Hence, for each node $i$, we obtain that
\begin{equation}
\min_{1\leq k\leq N}\left(  \lambda_{k}^{2}(A)\right)  =d_{i}\frac{1-S_{i}%
}{1+S_{i}}=d_{i}\left(  1-\frac{2}{\frac{1}{S_{i}}+1}\right)
\label{min_abs_eig_A}%
\end{equation}
Equation (\ref{min_abs_eig_A}) shows that $S_{i}\leq1$, and thus that $0\leq
S_{i}\leq1$ and that $\min_{1\leq k\leq N}\left(  \lambda_{k}^{2}(A)\right)
=d_{i}$ if $S_{i}=0$. Further, if $d_{i}>d_{j}$, then it follows from
(\ref{min_abs_eig_A}) that $S_{i}>S_{j}$ and $S_{i}$ increases with the degree
$d_{i}$. Hence, $S_{\min}=\min_{1\leq i\leq N}S_{i}$ corresponds to the node
with minimum degree.

Since each term in (\ref{def_Si}) is non-negative, $S_{i}$ can only be zero if
each term in the $k$-sum is zero,
\[
\left\{  \sum_{l=1}^{k}\left(  1-r_{i}\left(  l^{\ast}\right)  \right)
\right\}  \left(  \frac{1}{\frac{\lambda_{\left(  k+1\right)  }^{2}(A)}{d_{i}%
}+1}-\frac{1}{\frac{\lambda_{\left(  k\right)  }^{2}(A)}{d_{i}}+1}\right)  =0
\]
The first factor $\sum_{l=1}^{k}\left(  1-r_{i}\left(  l^{\ast}\right)
\right)  =\left(  1-r_{i}\left(  1\right)  \right)  +\sum_{l=2}^{k}\left(
1-r_{i}\left(  l^{\ast}\right)  \right)  $, because $r_{i}\left(  1^{\ast
}\right)  =r_{i}\left(  1\right)  $ as $\lambda_{\left(  1\right)  }%
^{2}(A)=\lambda_{1}^{2}(A)$. In a connected graph,
(\ref{i_th_component_xk_square}) demonstrates that $r_{i}\left(  1\right)
<1$, because each component of the principal eigenvector $x_{1}$ is positive
(by the Perron-Frobenius Theorem). Hence, for each $1\leq k\leq N-1$,
$\sum_{l=1}^{k}\left(  1-r_{i}\left(  l^{\ast}\right)  \right)  >0$. The last
factor cannot always be zero, because it would require that $\lambda_{\left(
k+1\right)  }^{2}(A)=\lambda_{\left(  k\right)  }^{2}(A)$ for all $k$, which
is impossible. Hence, in a connected graph, $S_{i}>0$ for each node $i$%
.\hfill$\square\medskip$

A consequence of Theorem \ref{theorem_bound_minimum_eigenvalue_A^2} is

\begin{corollary}
\label{corollary_disconnected_d_min_min_lambda_k} If $\min_{1\leq k\leq
N}\left(  \lambda_{k}^{2}(A)\right)  =d_{\min}$, then the graph is disconnected.
\end{corollary}

The reverse of the Corollary \ref{corollary_disconnected_d_min_min_lambda_k}
is not always true\footnote{Moreover, simulations on small Erd\H{o}s-R\'{e}nyi
graphs show that $\xi=d_{\min}-\min_{1\leq k\leq N}\left(  \lambda_{k}%
^{2}(A)\right)  -\mu_{N-1}$ is non-negative in most (but not all) cases.}.

If $\min_{1\leq k\leq N}\left(  \lambda_{k}^{2}(A)\right)  =0$, then $S_{i}=1$
for each node $i$. When excluding graphs with isolated nodes (i.e. degree zero
nodes), (\ref{min_abs_eig_A}) implies that, if $S_{i}=1$ for node $i$, then
$\min_{1\leq k\leq N}\left(  \lambda_{k}^{2}(A)\right)  =0$ and, thus
$S_{j}=1$ for each other node $j$. Hence, in any graph with $d_{i}>0$, in
order for $A$ to have a zero eigenvalue\footnote{It also follows from
(\ref{i_th_component_xk_square}) that, if $\lambda_{\left(  N\right)  }%
^{2}(A)=0$, then $(x_{N^{\ast}})_{i}^{2}=1-r_{i}\left(  N^{\ast}\right)  $,
where $N^{\ast}$ here equals the index $k$ for which $\lambda_{k}^{2}(A)=0$,
while summing over all $i$ shows that $\sum_{i=1}^{N}r_{i}\left(  N^{\ast
}\right)  =N-1$.}, there must hold that $S_{i}=1$ for each node $i$.

\section{The fundamental weight and its dual}

\label{sec_fundamental_weights} When choosing $b=u$ in Section
\ref{sec_determinantal_square_eigenvector_component}, the fundamental weight
$w_{k}=u^{T}x_{k}=\sum_{j=1}^{N}\left(  x_{k}\right)  _{j}$ was introduced as
additional information to determine the eigenvector components. The graph
angle $\gamma_{k}$ in \cite{Cvetkovic_eigenspaces} is related to the
fundamental weight by $\cos\gamma_{k}=\frac{w_{k}}{\sqrt{N}}$, where the angle
$\theta_{ab}$ between two vectors $a$ and $b$ obeys $\cos\theta_{ab}%
=\frac{a^{T}b}{\left\Vert a\right\Vert _{2}\left\Vert b\right\Vert _{2}}$ and
$\left\Vert a\right\Vert _{2}^{2}=a^{T}a$. Geometrically in $N=3$ dimensions,
the $3$ orthogonal axes are completely defined by the knowledge of 3 angles.
However, in higher dimensions $N>3$, all $\binom{N}{2}$ orthogonality
relations (\ref{first_orthogonality_equations}), which directly imply the
second set of $\binom{N}{2}$ orthogonality relations
(\ref{second_orthogonality_equations}) due commutativity between $X$ and
$X^{-1}$, are needed to specify the $N$ orthogonal axes ($x_{k}^{T}x_{m}%
=\cos\theta_{km}$ and $\theta_{km}$ is either $\frac{\pi}{2}$ or 0) so that we
expect $O\left(  N^{2}\right)  $ graph angles rather than $O\left(  N\right)
$.

Section \ref{sec_def_fundamental_weights} presents alternative definitions of
the fundamental and dual fundamental weights, while Section
\ref{sec_properties_fundamental_weights} derives a first set of their
properties. Using fundamental weights, we compute tight bounds on the coupling
of eigenvalues of a graph $G$ and its complement $G^{c}$ in Appendix
\ref{sec_eigenvectors_complementary_graph}.

\subsection{Definitions}

\label{sec_def_fundamental_weights}The dual of the definition
\begin{equation}
w_{k}=\sum_{j=1}^{N}\left(  x_{k}\right)  _{j}
\label{def_fundamental_weight_wk}%
\end{equation}
is%
\begin{equation}
\varphi_{j}=\sum_{k=1}^{N}\left(  x_{k}\right)  _{j}
\label{def_dual_fundamental_weight_phi_j}%
\end{equation}
which is the sum of the eigenvector components of node $j$ over all
eigenfrequencies. The corresponding vectors $w=\left(  w_{1},w_{2}%
,\ldots,w_{N}\right)  $ and $\varphi=\left(  \varphi_{1},\varphi_{2}%
,\ldots,\varphi_{N}\right)  $ are called the fundamental weight and dual
fundamental weight vector of the adjacency matrix $A$ of a graph $G$,
respectively. Those vectors can be written as the row sum and column sum of
the orthogonal matrix $X$ in (\ref{def_orthogonal_matrix_X}),
\begin{align}
w  &  =X^{T}u\label{def_fundamental_weight_vector_w}\\
\varphi &  =Xu \label{def_fundamental_weight_vector_phi}%
\end{align}
or, in terms of the eigenvectors $\left\{  x_{k}\right\}  _{1\leq k\leq N}$
and $\left\{  y_{k}\right\}  _{1\leq k\leq N}$, defined in
(\ref{def_y_nodal_eigenvector_components}),%
\begin{align}
w  &  =\sum_{k=1}^{N}y_{k}\label{w_sum_over_y}\\
\varphi &  =\sum_{k=1}^{N}x_{k} \label{phi_sum_over_eigenvectors_x}%
\end{align}
Hence, $\frac{1}{N}\varphi$ is the average of all eigenvectors of the
adjacency matrix $A$. The corresponding vector components are, for the
fundamental weight,%
\begin{equation}
w_{k}=u^{T}x_{k} \label{def_fundamental_weight_wk_scalar_product}%
\end{equation}
and for the dual fundamental weight%
\begin{equation}
\varphi_{j}=u^{T}y_{j}
\label{def_dual_fundamental_weight_phi_j_scalar_product}%
\end{equation}
illustrating that the role in (\ref{def_fundamental_weight_wk_scalar_product})
and (\ref{def_dual_fundamental_weight_phi_j_scalar_product}) of the vectors
$x_{k}$ and $y_{k}$ is reversed with respect to (\ref{w_sum_over_y}) and
(\ref{phi_sum_over_eigenvectors_x}).

Suppose that a node relabeling in the graph $G$ is defined by the permutation
matrix $P$, which is an orthogonal matrix obeying $Pu=u$. We denote the
relabeled adjacency matrix by $\widetilde{A}=P^{T}AP$ and its spectral
decomposition by $\widetilde{A}=\widetilde{X}\Lambda\widetilde{X}^{T}$, where
$\widetilde{X}=P^{T}X$. The definition (\ref{def_fundamental_weight_vector_w})
of $w$ shows that $\widetilde{w}=w$, so that $w$ is \emph{invariant} under a
relabeling transformation. However, the definition
(\ref{def_fundamental_weight_vector_phi}) of $\varphi$ shows that
$\widetilde{\varphi}=P^{T}\varphi$; in other words, the components of
$\widetilde{\varphi}$ change position after relabeling.

\begin{theorem}
\label{theorem_symmetric_X_exists.}There exist regular graphs for which the
adjacency matrix $A$ possesses a symmetric orthogonal matrix $X=X^{T}$.
\end{theorem}

Barik \emph{et al.} \cite{Barik_LAA2011} have shown that only regular graphs,
such as the complete graph $K_{N}$, for $N=4k$ and $k\in\mathbb{N}_{0}$, and
the regular bipartite graph $K_{2k,2k}$, are diagonalizable by a Hadamard
matrix. An $n\times n$ Hadamard matrix $H_{n}$ contains as elements either
$-1$ and $1$ and obeys $H_{n}H_{n}^{T}=nI_{n}$. The normalized matrix
$X_{n}=\frac{1}{\sqrt{n}}H_{n}$ is an orthogonal matrix, from which it follows
that $\det H_{n}=n^{\frac{n}{2}}$, which is maximal among all $n\times n$
matrices with elements in absolute value less than or equal to $1$, which
includes all orthogonal matrices. Any relabeling (permutation of rows and
columns) of a Hadamard matrix is again a Hadamard matrix; multiplying any row
or column by $-1$ preserves the Hadamard properties.

Sylvester found a construction for \emph{symmetric} Hadamard matrices
$H_{2^{k}}=H_{2^{k-1}}\otimes H_{2}$, where $\otimes$ is the Kronecker product
\cite{PVM_graphspectra} and $H_{2}=\left[
\begin{array}
[c]{cc}%
1 & 1\\
1 & -1
\end{array}
\right]  $, that contain the $u$ vector in the first column.

\textbf{Proof\footnote{The proof (only for $K_{N}$) is slightly simpler than
the one in \cite{Barik_LAA2011}.}:} Let $H_{n}=\left[  u|\widetilde{H}\right]
$ so that $H_{n}e_{1}=u$. Consider the diagonal matrix $D=I-e_{1}e_{1}^{T}$,
then%
\[
H_{n}DH_{n}^{T}=H_{n}H_{n}^{T}-H_{n}e_{1}\left(  H_{n}e_{1}\right)
^{T}=nI_{n}-u.u^{T}=nI-J
\]
Hence, the Laplacian matrix of the complete graph $K_{n}$ is $Q_{K_{n}%
}=nI-J=H_{n}DH_{n}^{T}$. Since $K_{n}$ is a regular graph, the eigenvectors of
the Laplacian $Q$ and the adjacency matrix $A$ are the same\footnote{Indeed,
for a regular graph with degree $r$, the Laplacian is $Q=\left(  r+1\right)
I-A$. If $Q=ZMZ^{T}$ and $A=X\Lambda X^{T}$, we observe that $ZMZ^{T}=X\left(
\left(  r+1\right)  I-\Lambda\right)  X^{T}$, implying that $X=Z$.}. In
conclusion, any Hadamard matrix with $H_{n}e_{1}=u$ provides the orthogonal
matrix for the complete graph (where $N=4k$) and the Sylvester construction
demonstrates that there exist symmetric such Hadamard matrices.\hfill
$\square\medskip$

It will transpire that graphs with a symmetric orthogonal matrix $X=X^{T}$
possess extremal properties: (\ref{def_fundamental_weight_vector_w}) and
(\ref{def_fundamental_weight_vector_w}) show that $w=$ $\varphi$, only if
$X=X^{T}$.

\subsection{Properties}

\label{sec_properties_fundamental_weights}The definitions in Section
\ref{sec_def_fundamental_weights} lead to a number of immediate consequences.

\medskip\medskip\refstepcounter{article}{\noindent\textbf{\thearticle. }%
}\ignorespaces\label{art_norm_w_phi} The norms $\left\Vert w\right\Vert _{2}$
and $\left\Vert \varphi\right\Vert _{2}$ are the same. In particular,%
\[
w^{T}w=\varphi^{T}\varphi=N
\]
shows that their norm equals that of the all-one vector $u$,%
\[
\left\Vert w\right\Vert _{2}=\left\Vert \varphi\right\Vert _{2}=\left\Vert
u\right\Vert _{2}=\sqrt{N}%
\]
This norm property follows either from (\ref{def_fundamental_weight_vector_w})
as $w^{T}w=u^{T}XX^{T}u=u^{T}u$, because $XX^{T}=I$, or from $w^{T}%
w=\sum_{k=1}^{N}\sum_{m=1}^{N}y_{k}^{T}y_{m}$ and the orthogonality relations
(\ref{second_orthogonality_equations}). Similarly for $\varphi$, where a
possible node relabeling does not influence the norm: $\widetilde{\varphi}%
^{T}\widetilde{\varphi}=\varphi^{T}PP^{T}\varphi=\varphi^{T}\varphi$, because
$P$ is an orthogonal matrix.

\medskip\refstepcounter{article}{\noindent\textbf{\thearticle. }%
}\ignorespaces\label{art_elementary_bound}\emph{Bounds of maximum and
minimum}.\textbf{ }The bound $N=\varphi^{T}\varphi=\sum_{j=1}^{N}\varphi
_{j}^{2}\leq N\max_{1\leq j\leq N}\varphi_{j}^{2}$ illustrates that%
\[
-\sqrt{N}\leq\min_{1\leq j\leq N}\varphi_{j}\leq1\leq\max_{1\leq j\leq
N}\varphi_{j}\leq\sqrt{N}%
\]
and similarly for $w$. In the case of $w$, a much sharper lower bound for the
maximum is known (\textbf{art. }\ref{art_bounds_phi_w}).

\medskip\refstepcounter{article}{\noindent\textbf{\thearticle. }%
}\ignorespaces\label{art_projection_w_on_ym_phi_on_xm}\textbf{ }Any vector $z$
in an $N$-dimensional space can be written as a linear combination of a set of
$N$ orthogonal vectors that span that space, such as the set $\left\{
x_{k}\right\}  _{1\leq k\leq N}$ and the set $\left\{  y_{k}\right\}  _{1\leq
k\leq N}$,%
\[
z=\sum_{k=1}^{N}\left(  z^{T}y_{k}\right)  y_{k}=\sum_{k=1}^{N}\left(
z^{T}x_{k}\right)  x_{k}%
\]
For example,
\begin{equation}
u=\sum_{j=1}^{N}\varphi_{j}y_{j}=\sum_{k=1}^{N}w_{k}x_{k}
\label{u_as_lin_comb_eigenvectors}%
\end{equation}
indicating that the coordinate vector of $u$ with respect to the basis
$\left\{  x_{k}\right\}  _{1\leq k\leq N}$ is the $w$ vector and with respect
to the basis $\left\{  y_{j}\right\}  _{1\leq j\leq N}$ is the $\varphi$
vector. Another example is, using (\ref{third_orthogonality_combination}),%
\[
x_{k}=\sum_{j=1}^{N}\left(  x_{k}^{T}y_{j}\right)  y_{j}=\sum_{j=1}^{N}\left(
X^{2}\right)  _{jk}y_{j}\text{ and }y_{j}=\sum_{k=1}^{N}\left(  X^{2}\right)
_{jk}x_{k}%
\]

The definitions (\ref{w_sum_over_y}) and (\ref{phi_sum_over_eigenvectors_x})
express $w$ and $\varphi$ as such a linear combination, from which it follows,
for any integer $1\leq m\leq N$, that%
\begin{equation}
w^{T}y_{m}=1 \label{scalar_product_w_y_k=1}%
\end{equation}
and\footnote{We also mention the dual expressions, derived by invoking
(\ref{third_orthogonality_combination}),%
\[
\left\{
\begin{array}
[c]{c}%
w^{T}x_{m}=\left(  X^{2}u\right)  _{m}\\
\varphi^{T}y_{m}=\left(  X^{2}u\right)  _{m}%
\end{array}
\right.
\]
}%
\begin{equation}
\varphi^{T}x_{m}=1 \label{scalar_product_phi_xk=1}%
\end{equation}
Both scalar products (\ref{scalar_product_w_y_k=1}) and
(\ref{scalar_product_phi_xk=1}) also follow from the identities $XX^{T}u=u$
and $X^{T}Xu=u$, respectively. The geometric meaning is that, for any $m$, the
vector $w$ and $\varphi$ make the same angle $\theta_{y_{m}w}$ and
$\theta_{x_{m}\varphi}$ with any vector $y_{m}$ and any eigenvector $x_{m}$ of
$A$, respectively. Hence, with respect to the orthogonal basis spanned by the
eigenvectors $\left\{  x_{k}\right\}  _{1\leq k\leq N}$, the vector $\varphi$
plays the same role as the vector $u$ with respect to the \textquotedblleft
classical\textquotedblright\ orthogonal basis $\left\{  e_{k}\right\}  _{1\leq
k\leq N}$. The transformation $Xu=\varphi$ rotates the all-one vector $u$ into
the vector $\varphi$, while the inverse rotation yields $X^{-1}u=X^{T}u=w$.

\medskip\refstepcounter{article}{\noindent\textbf{\thearticle. }%
}\ignorespaces\label{art_lambda^T.w=d^T.phi}\textbf{ }After left-multiplying
the eigenvalue equation $Ax_{k}=\lambda_{k}x_{k}$ by $u^{T}$ and summing the
resulting eigenvalue relation $d^{T}x_{k}=\lambda_{k}w_{k}$ over all $k$
yields%
\[
\sum_{k=1}^{N}\lambda_{k}w_{k}=\sum_{k=1}^{N}\left(  \sum_{j=1}^{N}%
d_{j}\left(  x_{k}\right)  _{j}\right)  =\sum_{j=1}^{N}d_{j}\varphi_{j}%
\]
In other words, we observe that
\begin{equation}
w^{T}\lambda=\varphi^{T}d \label{scalar_prod_w_lambda=d_phi}%
\end{equation}
where the eigenvalue vector $\lambda=\left(  \lambda_{1},\lambda_{2}%
,\ldots,\lambda_{N}\right)  $ is related to the degree vector $d$ via the
fundamental weight vector $w$ and its dual vector $\varphi$. Further, recall
that $\lambda^{T}\lambda=2L$ and $d^{T}u=2L$, but $d^{T}d>2L$. Thus,
$w^{T}\lambda=\varphi^{T}d$, combined with%
\[
\left(  d^{T}\varphi\right)  ^{2}\leq\left\Vert d\right\Vert _{2}%
^{2}\left\Vert \varphi\right\Vert _{2}^{2}=\left\Vert d\right\Vert _{2}^{2}N
\]
and%
\[
\left(  \lambda^{T}w\right)  ^{2}\leq\left\Vert \lambda\right\Vert _{2}%
^{2}\left\Vert w\right\Vert _{2}^{2}=2LN\leq\left\Vert d\right\Vert _{2}^{2}N
\]
means that the angle $\theta_{d\varphi}$ between the vector $d$ and $\varphi$
is larger than the angle $\theta_{\lambda w}$ between the vector $\lambda$ and
$w$. Thus, we can say that $\lambda$ and $w$ are closer correlated than $d$
and $\varphi$.

The generalization of (\ref{scalar_prod_w_lambda=d_phi}), based on the
eigenvalue equation $A^{m}x_{k}=\lambda_{k}^{m}x_{k}$ is%
\begin{equation}
w^{T}\lambda^{m}=\varphi^{T}d^{\left(  m\right)  }
\label{scalar_prod_w_lambda=d_phi_m}%
\end{equation}
or%
\[
\sum_{k=1}^{N}\lambda_{k}^{m}w_{k}=\sum_{j=1}^{N}\left(  A^{m}u\right)
_{j}\varphi_{j}%
\]
where the vector $\lambda^{m}=$ $\left(  \lambda_{1}^{m},\lambda_{2}%
^{m},\ldots,\lambda_{N}^{m}\right)  =\Lambda^{m}u$ and the vector $d^{\left(
m\right)  }=\left(  \left(  A^{m}u\right)  _{1},\left(  A^{m}u\right)
_{2},\ldots,\left(  A^{m}u\right)  _{N}\right)  =A^{m}u$, with $d^{\left(
0\right)  }=u$ and $d^{\left(  1\right)  }=d$.

Although the vector $d$ cannot be equal to the vector $\lambda$, we cannot
conclude from (\ref{scalar_prod_w_lambda=d_phi_m}) that $w$ cannot be equal to
$\varphi$. Indeed, suppose that $w=\varphi\neq0$, then
(\ref{scalar_prod_w_lambda=d_phi_m}) reduces to $w^{T}\left(  \lambda
^{m}-d^{\left(  m\right)  }\right)  =0$, which is equivalent to
\[
0=w^{T}\left(  A^{m}-\Lambda^{m}\right)  u=w^{T}\left(  X\Lambda^{m}%
X^{T}-\Lambda^{m}XX^{T}\right)  u=w^{T}\left(  X\Lambda^{m}-\Lambda
^{m}X\right)  w=-w^{T}\left(  X^{T}\Lambda^{m}-\Lambda^{m}X^{T}\right)  w
\]
and only possible for all $m$ if $X=X^{T}$.

\medskip\refstepcounter{article}{\noindent\textbf{\thearticle. }%
}\ignorespaces\label{art_sum_elem_X}\textbf{ } Since $s_{X}=u^{T}Xu=\left(
X^{T}u\right)  ^{T}u=u^{T}\left(  Xu\right)  $, we have with
(\ref{def_fundamental_weight_vector_w}) and
(\ref{def_fundamental_weight_vector_phi}) that%
\begin{equation}
s_{X}=w^{T}u=u^{T}\varphi\label{sum_elements_orthogonal_matrix_X}%
\end{equation}
which also follows from (\ref{scalar_prod_w_lambda=d_phi_m}) for $m=0$. The
sum of the elements of the orthogonal matrix $X$ thus equals
\[
s_{X}=\sum_{j=1}^{N}\sum_{k=1}^{N}\left(  x_{k}\right)  _{j}=\sum_{j=1}%
^{N}\varphi_{j}=\sum_{k=1}^{N}w_{k}%
\]
Since $\left\vert s_{X}\right\vert =\left\vert w^{T}u\right\vert \leq
\sqrt{\left\Vert w\right\Vert _{2}^{2}\left\Vert u\right\Vert _{2}^{2}}$, we
find with $\left\Vert w\right\Vert _{2}^{2}=N$ that $-N\leq s_{X}\leq N$.
Moreover, the sum of the elements of the matrix $ww^{T}=X^{T}JX$ and its
transpose $\varphi\varphi^{T}=XJX^{T}$ (where the all-one matrix is
$J=u.u^{T}$) equals
\[
s_{X}^{2}=u^{T}\left(  ww^{T}\right)  u=u^{T}\left(  \varphi\varphi
^{T}\right)  u
\]

\medskip\refstepcounter{article}{\noindent\textbf{\thearticle. }%
}\ignorespaces\label{art_sum_elem_X^2}\textbf{ }Since $X^{T}=X^{-1}$ is
non-singular ($\det X=\pm1$), it follows from
(\ref{def_fundamental_weight_vector_w}) and
(\ref{def_fundamental_weight_vector_phi}) that the all-one vector can be
expressed as
\[
u=Xw=X^{T}\varphi
\]
Thus, we find from the definition (\ref{def_fundamental_weight_vector_w}) and
(\ref{def_fundamental_weight_vector_phi}) that%
\[
\varphi=X^{2}w\text{ and }w=\left(  X^{T}\right)  ^{2}\varphi
\]
and
\[
w^{T}\varphi=w^{T}X^{2}w=u^{T}X^{2}u
\]
so that the sum $s_{X^{2}}=u^{T}X^{2}u$ of the elements in $X^{2}$ equals%
\begin{equation}
s_{X^{2}}=w^{T}\varphi\label{sum_elements_X^2}%
\end{equation}
Since $w^{T}\varphi=\left\Vert w\right\Vert _{2}\left\Vert \varphi\right\Vert
_{2}\cos\theta_{w\varphi}=N\cos\theta_{w\varphi}$, we find that the sum of the
elements of $X^{2}$ obeys $-N\leq u^{T}X^{2}u\leq N$ and equality in either
lower or upper bound only holds if $w=-\varphi$ or $w=\varphi$.

\medskip\refstepcounter{article}{\noindent\textbf{\thearticle. }%
}\ignorespaces\label{art_closed_walks_and_phi}\textbf{ }Since $u^{T}\lambda=0$
due to trace$\left(  A\right)  =0$, we find
\begin{equation}
\varphi^{T}A\varphi=0 \label{phi_A_phi = 0}%
\end{equation}
Relation (\ref{phi_A_phi = 0}) holds for any adjacency matrix $A=A^{T}$. At
first glance, the Rayleigh equations may hint that $\varphi$ is an eigenvector
belonging to eigenvalue $\lambda=0$, which is false, because for any component
$j$, we find
\[
\left(  A\varphi\right)  _{j}=\left(  AXu\right)  _{j}=\left(  X\Lambda
u\right)  _{j}=\left(  X\lambda\right)  _{j}=\sum_{k=1}^{N}\left(
x_{k}\right)  _{j}\lambda_{k}%
\]
while (\ref{zero_node_j_square_eigenvector_components}) shows that only
$\sum_{k=1}^{N}\left(  x_{k}\right)  _{j}^{2}\lambda_{k}=0$. Moreover,
(\ref{phi_A_phi = 0}) demonstrates that not all components of $\varphi$ can be
negative nor all can be positive for a graph (except for the null graph
without any links for which $A=0$). For, otherwise, if $\varphi_{j}<0$ or
$\varphi_{j}>0$ for all $j$, then $\varphi^{T}A\varphi=2\sum_{i=1}^{N}%
\sum_{j=i+1}^{N}a_{ij}\varphi_{i}\varphi_{j}>0$ contradicting
(\ref{phi_A_phi = 0}). The vector $\varphi=\sqrt{N}e_{j}$ (see \textbf{art.
}\ref{art_w_phi_regular_graphs}) does not violate (\ref{phi_A_phi = 0})
because $a_{jj}=0$. In general, $\varphi$ has positive, zero and negative
components. It is convenient to order (e.g. by a node relabeling) the dual
fundamental weights as
\[
\varphi_{\left(  1\right)  }\geq\ldots\geq\varphi_{\left(  q\right)  }%
\geq0\geq\varphi_{\left(  q+1\right)  }\geq\ldots\geq\varphi_{\left(
N\right)  }%
\]
with $1\leq q<N$.

A generalization of (\ref{phi_A_phi = 0}) follows from $A^{m}=X\Lambda
^{m}X^{T}$ and $u=X^{T}\varphi$ as
\[
\varphi^{T}A^{m}\varphi=\varphi^{T}X\Lambda^{m}X^{T}\varphi=u^{T}\Lambda
^{m}u=\sum_{j=1}^{N}\lambda_{j}^{m}=W_{m}%
\]
Hence, the number $W_{m}$ of closed walks with $m$ hops equals%
\begin{equation}
W_{m}=\varphi^{T}A^{m}\varphi=\text{ trace}\left(  A^{m}\right)  =\sum
_{j=1}^{N}\lambda_{j}^{m} \label{number_closed_k_hops_walks_phi}%
\end{equation}
whereas the total number $N_{m}$ of walks with $m$ hops equals%
\begin{equation}
N_{m}=u^{T}A^{m}u=w^{T}\Lambda^{m}w=\sum_{j=1}^{N}w_{j}^{2}\lambda_{j}^{m}
\label{number_all_m_hops_walks}%
\end{equation}

\medskip\refstepcounter{article}{\noindent\textbf{\thearticle. }%
}\ignorespaces\label{art_w_phi_regular_graphs} \emph{Regular graphs.} In any
connected regular graph (with degree vector $d=ru$), it holds that
$x_{1}=\frac{u}{\sqrt{N}}$. Since $u$ is always an eigenvector of the
Laplacian matrix $Q$, this \textbf{art. }\ref{art_w_phi_regular_graphs} holds
for the Laplacian of any graph as well (with replacement below of $e_{1}$ by
$e_{N}$, because $u$ corresponds to the smallest Laplacian eigenvalue $\mu
_{N}=0$). If the graph is not connected, a different normalization of $u$ is required.

The definition (\ref{def_fundamental_weight_wk_scalar_product}) indicates that
$w_{1}=u^{T}x_{1}=\sqrt{N}$, while $w_{k}=0$ (due to orthogonality
(\ref{first_orthogonality_equations})), so that the entire $w=\sqrt{N}e_{1}$
vector is known. Thus, if the graph is regular, then $w=\sqrt{N}e_{1}$ and the
sum (\ref{sum_elements_orthogonal_matrix_X}) of the elements of $X$ equals
$s_{X}=u^{T}w=\sqrt{N}$, while (\ref{sum_elements_X^2}) shows that $s_{X^{2}%
}=$ $\varphi_{1}$, which is clearly not invariant to node relabeling! The
converse, \textquotedblleft if $s_{X}=\sqrt{N}$, then the graph is
regular\textquotedblright, is likely not true\footnote{A counter example of an
irregular graph with $s_{X}\simeq$ $\sqrt{N}$ (up to 6 digits accurate) has
been found by Xiangrong Wang.}.

\medskip\refstepcounter{article}{\noindent\textbf{\thearticle. }%
}\ignorespaces\label{art_bounds_phi_w}\emph{Bounds.} The definition
$w_{k}=u^{T}x_{k}=\sum_{j=1}^{N}\left(  x_{k}\right)  _{j}$ and the first
orthogonality condition (\ref{first_orthogonality_equations}) yield%
\[
w_{k}=\frac{u^{T}x_{k}}{x_{k}^{T}x_{k}}=\frac{\sum_{j^{\prime}=1}^{N}\left(
x_{k}\right)  _{j^{\prime}}}{\sum_{j^{\prime}=1}^{N}\left(  x_{k}\right)
_{j^{\prime}}^{2}}%
\]
where $j^{\prime}$ reflects that $\left(  x_{k}\right)  _{j^{\prime}}\neq0$,
because a zero term does not contribute to the sum. The inequality
\cite{Hardy_inequality}
\begin{equation}
\min_{1\leq j\leq n}\frac{r_{j}}{b_{j}}\leq\frac{r_{1}+r_{2}+\cdots+r_{n}%
}{b_{1}+b_{2}+\cdots+b_{n}}\leq\max_{1\leq j\leq n}\frac{r_{j}}{b_{j}}
\label{inequality_breuk_n_termen}%
\end{equation}
where $b_{1},b_{2},\ldots,b_{n}$ are positive real numbers and $r_{1}%
,r_{2},\ldots,r_{n}$ are real numbers, yields%
\begin{equation}
\min_{1\leq j^{\prime}\leq N}\frac{1}{\left(  x_{k}\right)  _{j^{\prime}}}\leq
w_{k}\leq\max_{1\leq j^{\prime}\leq N}\frac{1}{\left(  x_{k}\right)
_{j^{\prime}}} \label{inequality_fundamental_weight}%
\end{equation}
and, similarly,
\[
\min_{1\leq k^{\prime}\leq N}\frac{1}{\left(  x_{k^{\prime}}\right)  _{j}}%
\leq\varphi_{j}\leq\max_{1\leq k^{\prime}\leq N}\frac{1}{\left(  x_{k^{\prime
}}\right)  _{j}}%
\]

All components of $x_{1}$ are non-negative by the Perron-Frobenius theorem,
whereas $\min_{1\leq j^{\prime}\leq N}\frac{1}{\left(  x_{k}\right)
_{j^{\prime}}}<-1$ for $k>1$, because any other eigenvector $x_{k}$ must have
at least one negative component to satisfy the orthogonality condition
$x_{k}^{T}x_{1}=0$. The inequality\ (\ref{inequality_fundamental_weight})
illustrates that $w_{1}\geq1$, a result earlier found in \cite[p.
40]{PVM_graphspectra} with a different method. A much sharper lower bound%
\begin{equation}
\sqrt{\frac{\lambda_{1}}{1-\frac{1}{\omega}}}\leq w_{1}
\label{bounds_fundamental_weight_w1}%
\end{equation}
where $\omega$ is the clique number of the graph $G$ is proved and evaluated
in \cite{PVM_Lower_bound_fundamental_weight_SITIS2014}.

\medskip\refstepcounter{article}{\noindent\textbf{\thearticle. }%
}\ignorespaces\label{art_minimal_spacing_phi}\textbf{ }\emph{Upper bound for
the minimal spacing}. We assume that the vector components of $\varphi$ are
ordered as in \textbf{art. }\ref{art_closed_walks_and_phi}. The corresponding
relabeling\footnote{Notice that $d_{\left(  l\right)  }$ denotes the $l$-th
largest degree in the graph, while $d_{l^{\ast}}$ is the degree of the node
$l^{\ast}$, whose dual fundamental weight component $\varphi_{l^{\ast}%
}=\varphi_{\left(  l\right)  }$ is the \thinspace$l$-th largest.} of node $l$
is denoted by $l^{\ast}$. Several bounds for the spacings $\varphi_{\left(
j\right)  }-\varphi_{\left(  j+1\right)  }$ will be derived, based on Theorem
\ref{theo_spacing_sequence} in Appendix \ref{sec_spacings}.

We apply (\ref{def_f}) to $c=\varphi-u\varphi_{\left(  N\right)  }$ with
$a=e_{1}$ (telescoping series) and with $a=u$ so that $a^{T}c=s_{X}%
-N\varphi_{\left(  N\right)  }$ in (\ref{sum_elements_orthogonal_matrix_X}).
The corresponding fractions $f$ are%
\begin{align*}
f_{e_{1}}  &  =\frac{\varphi_{\left(  1\right)  }-\varphi_{\left(  N\right)
}}{N-1}\\
f_{u}  &  =\frac{\frac{s_{X}}{N}-\varphi_{\left(  N\right)  }}{\frac{N-1}{2}}%
\end{align*}
Since $\left\vert \varphi_{\left(  1\right)  }\right\vert \leq\sqrt{N}$ and
$\left\vert \varphi_{\left(  N\right)  }\right\vert \leq\sqrt{N}$
(\textbf{art. }\ref{art_elementary_bound}), we conclude from Theorem
\ref{theo_spacing_sequence} that
\[
\min_{1\leq j\leq N-1}\left\{  \varphi_{\left(  j\right)  }-\varphi_{\left(
j+1\right)  }\right\}  \leq O\left(  \frac{1}{\sqrt{N}}\right)
\]
for large $N$.

\medskip\refstepcounter{article}{\noindent\textbf{\thearticle. }%
}\ignorespaces\label{art_maximal_spacing_phi}\textbf{ }\emph{Lower bound for
the maximal spacing}. For $a=\varphi$ and $\varphi^{T}\varphi=N$, the fraction
(\ref{def_f}) becomes
\[
f_{\varphi}=\frac{N-\varphi_{\left(  N\right)  }s_{X}}{\sum_{m=1}%
^{N-1}m\varphi_{\left(  N-m\right)  }}%
\]
Using the Cauchy-Schwarz inequality \cite[p. 257]{PVM_graphspectra}%
\[
\sum_{m=1}^{N-1}m\varphi_{\left(  N-m\right)  }\leq\sqrt{\sum_{m=1}^{N-1}%
m^{2}\sum_{m=1}^{N-1}\varphi_{\left(  N-m\right)  }^{2}}=\sqrt{\frac{\left(
N-1\right)  N\left(  2N-1\right)  }{6}\left(  N-\varphi_{\left(  N\right)
}^{2}\right)  }%
\]
we arrive at%
\[
\max_{1\leq j\leq N-1}\left\{  \varphi_{\left(  j\right)  }-\varphi_{\left(
j+1\right)  }\right\}  \geq\frac{1-\varphi_{\left(  N\right)  }\frac{s_{X}}%
{N}}{\sqrt{\frac{\left(  N-1\right)  \left(  2N-1\right)  }{6}\left(
1-\frac{\varphi_{\left(  N\right)  }^{2}}{N}\right)  }}\geq O\left(  \frac
{1}{N}\right)
\]
We observe that sharper bounds here and in \textbf{art. }%
\ref{art_minimal_spacing_phi} are only possible when $\varphi_{\left(
N\right)  }$ is known. After applying Abel summation to
(\ref{sum_elements_orthogonal_matrix_X}),%
\[
\sum_{j=1}^{N-1}j\left\{  \varphi_{\left(  j\right)  }-\varphi_{\left(
j+1\right)  }\right\}  =s_{X}-N\varphi_{\left(  N\right)  }%
\]
we find (since $\varphi_{\left(  j\right)  }-\varphi_{\left(  j+1\right)
}\geq0$) that $\frac{s_{X}}{N}\geq\varphi_{\left(  N\right)  }\geq-\sqrt{N}$,
where $s_{X}$ can be negative.

\medskip\refstepcounter{article}{\noindent\textbf{\thearticle. }%
}\ignorespaces\label{art_maximal_spacing_phi_another_bound} \emph{Another type
of lower bound for the maximal spacing}. \textbf{Art. }%
\ref{art_lambda^T.w=d^T.phi} demonstrates that $d^{T}\varphi=\lambda^{T}w$ and
$-\sqrt{2LN}\leq\lambda^{T}w\leq\sqrt{2LN}$, so that, after Abel summation,%
\[
-\sqrt{2LN}\leq\sum_{j=1}^{N-1}\left(  \sum_{l=1}^{j}d_{l^{\ast}}\right)
\left\{  \varphi_{\left(  j\right)  }-\varphi_{\left(  j+1\right)  }\right\}
+2L\varphi_{\left(  N\right)  }\leq\sqrt{2LN}%
\]
or%
\[
-\sqrt{\frac{1}{d_{av}}}-\frac{1}{2L}\sum_{j=1}^{N-1}\left(  \sum_{l=1}%
^{j}d_{l^{\ast}}\right)  \left\{  \varphi_{\left(  j\right)  }-\varphi
_{\left(  j+1\right)  }\right\}  \leq\varphi_{\left(  N\right)  }\leq
\sqrt{\frac{1}{d_{av}}}-\frac{1}{2L}\sum_{j=1}^{N-1}\left(  \sum_{l=1}%
^{j}d_{l^{\ast}}\right)  \left\{  \varphi_{\left(  j\right)  }-\varphi
_{\left(  j+1\right)  }\right\}
\]
Since $\varphi_{\left(  N\right)  }<0$ (ignoring pathological cases), we find
that%
\[
\sqrt{\frac{1}{d_{av}}}<\frac{1}{2L}\sum_{j=1}^{N-1}\left(  \sum_{l=1}%
^{j}d_{l^{\ast}}\right)  \left\{  \varphi_{\left(  j\right)  }-\varphi
_{\left(  j+1\right)  }\right\}
\]
Further,%
\[
\sum_{j=1}^{N-1}\left(  \sum_{l=1}^{j}d_{l^{\ast}}\right)  \left\{
\varphi_{\left(  j\right)  }-\varphi_{\left(  j+1\right)  }\right\}  \leq
\max_{1\leq j\leq N-1}\left\{  \varphi_{\left(  j\right)  }-\varphi_{\left(
j+1\right)  }\right\}  \sum_{j=1}^{N-1}\left(  \sum_{l=1}^{j}d_{l^{\ast}%
}\right)
\]
and%
\[
\sum_{j=1}^{N-1}\sum_{l=1}^{j}d_{l^{\ast}}=\sum_{l=1}^{N-1}\sum_{j=l}%
^{N-1}d_{l^{\ast}}=\sum_{l=1}^{N-1}\left(  N-l^{\ast}\right)  d_{l^{\ast}%
}=\sum_{l=1}^{N}\left(  N-l^{\ast}\right)  d_{l^{\ast}}=2LN-\sum_{l=1}%
^{N}l^{\ast}d_{l^{\ast}}%
\]
Hence,%
\[
\max_{1\leq j\leq N-1}\left\{  \varphi_{\left(  j\right)  }-\varphi_{\left(
j+1\right)  }\right\}  >\frac{\sqrt{\frac{1}{d_{av}}}}{N-\frac{1}{2L}%
\sum_{l=1}^{N}l^{\ast}d_{l^{\ast}}}>\frac{\sqrt{\frac{1}{d_{av}}}}%
{N-\frac{d_{\min}}{d_{av}}\frac{\left(  N+1\right)  }{2}}>\frac{1}{N}%
\frac{\sqrt{\frac{1}{d_{av}}}}{1-\frac{d_{\min}}{2d_{av}}}%
\]
which is, unfortunately, more conservative than the lower bound in
\textbf{art. }\ref{art_maximal_spacing_phi}.

\section{Conclusion}

\label{sec_conclusion}Three Theorems
\ref{theorem_eigenvector_component_determinants},
\ref{theorem_square_eigenvector_component_error_r} and
\ref{theorem_eigenvectors_in_terms_walks} present different expressions for
the square of eigenvector components of the adjacency matrix of a graph. Many
other formulae and bounds are deduced from those Theorems. Section
\ref{sec_interpretations} proposes the fundamental expression
(\ref{squared_eigenvector_component_xk_j_charc_pol}) as a nodal centrality
metric and shows its relation to the notion of graph robustness. Section
\ref{sec_fundamental_weights} presents the definition and properties of the
fundamental weights and the dual fundamental weights.\medskip

\textbf{Acknowledgements}. I am very grateful to Willem Haemers, Edwin van Dam
and Dragos Cvetcovic for their input.

{\footnotesize
\bibliographystyle{unsrt}
\bibliography{cac,MATH,misc,net,pvm,QTH,tel}
}

\appendix{}

\section{Eigenvectors and eigenvalues: brief review}

\label{sec_introduction_eigenvectors_eigenvalues}

\subsection{Definition}

We denote by $x_{k}$ the eigenvector of the symmetric matrix $A$ belonging to
the eigenvalue $\lambda_{k}$, normalized so that $x_{k}^{T}x_{k}=1$. The
eigenvalues of an $N\times N$ symmetric matrix $A=A^{T}$ are real and can be
ordered as $\lambda_{1}\geq\lambda_{2}\geq\ldots\geq\lambda_{N}$. Let $X$ be
the orthogonal matrix with eigenvectors of $A$ in the columns,%
\[
X=\left[
\begin{array}
[c]{ccccc}%
x_{1} & x_{2} & x_{3} & \cdots & x_{N}%
\end{array}
\right]
\]
or explicitly in terms of the $m$-th component $\left(  x_{j}\right)  _{m}$ of
eigenvector $x_{j}$,%
\begin{equation}
X=\left[
\begin{array}
[c]{ccccc}%
\left(  x_{1}\right)  _{1} & \left(  x_{2}\right)  _{1} & \left(
x_{3}\right)  _{1} & \cdots & \left(  x_{N}\right)  _{1}\\
\left(  x_{1}\right)  _{2} & \left(  x_{2}\right)  _{2} & \left(
x_{3}\right)  _{2} & \cdots & \left(  x_{N}\right)  _{2}\\
\left(  x_{1}\right)  _{3} & \left(  x_{2}\right)  _{3} & \left(
x_{3}\right)  _{3} & \cdots & \left(  x_{N}\right)  _{3}\\
\vdots & \vdots & \vdots & \ddots & \vdots\\
\left(  x_{1}\right)  _{N} & \left(  x_{2}\right)  _{N} & \left(
x_{3}\right)  _{N} & \cdots & \left(  x_{N}\right)  _{N}%
\end{array}
\right]  \label{def_orthogonal_matrix_X}%
\end{equation}
where the element $X_{ij}=\left(  x_{j}\right)  _{i}$. The eigenvalue equation
$Ax_{k}=\lambda_{k}x_{k}$ translates to the matrix equation $A=X\Lambda X^{T}%
$, where $\Lambda=$ diag$\left(  \lambda_{k}\right)  $.

The relation $X^{T}X=I=XX^{T}$ (see e.g. \cite[p. 223]{PVM_graphspectra})
expresses, in fact, \emph{double orthogonality}. The first equality $X^{T}X=I$
translates to the well-known orthogonality relation
\begin{equation}
x_{k}^{T}x_{m}=\sum_{j=1}^{N}\left(  x_{k}\right)  _{j}\left(  x_{m}\right)
_{j}=\delta_{km} \label{first_orthogonality_equations}%
\end{equation}
stating that the eigenvector $x_{k}$ belonging to eigenvalue $\lambda_{k}$ is
orthogonal to any other eigenvector belonging to a different eigenvalue. The
second equality $XX^{T}=I$, which arises from the commutativity of the inverse
matrix $X^{-1}=X^{T}$ with the matrix $X$ itself, can be written as
$\sum_{j=1}^{N}\left(  x_{j}\right)  _{m}\left(  x_{j}\right)  _{k}%
=\delta_{mk}$ and suggests us to define the row vector in $X$ as
\begin{equation}
y_{m}=\left(  \left(  x_{1}\right)  _{m},\left(  x_{2}\right)  _{m}%
,\ldots,\left(  x_{N}\right)  _{m}\right)
\label{def_y_nodal_eigenvector_components}%
\end{equation}
Then, the second orthogonality condition $XX^{T}=I$ implies orthogonality of
the vectors%
\begin{equation}
y_{l}^{T}y_{j}=\sum_{k=1}^{N}(x_{k})_{l}(x_{k})_{j}=\delta_{lj}
\label{second_orthogonality_equations}%
\end{equation}
Beside the first (\ref{first_orthogonality_equations}) and second
(\ref{second_orthogonality_equations}) orthogonality relations, the third
combination equals
\begin{equation}
y_{j}^{T}x_{k}=\sum_{l=1}^{N}\left(  x_{l}\right)  _{j}\left(  x_{k}\right)
_{l}=\sum_{l=1}^{N}X_{jl}X_{lk}=\left(  X^{2}\right)  _{jk}
\label{third_orthogonality_combination}%
\end{equation}

\subsection{Frequency interpretation}

The sum over $j$ in (\ref{second_orthogonality_equations}) can be interpreted
as the sum over all eigenvalues. Indeed, the eigenvalue equation is%
\begin{equation}
Ax\left(  \lambda\right)  =\lambda\ x\left(  \lambda\right)
\label{eigenvalue_equation_lambda}%
\end{equation}
where a non-zero vector $x\left(  \lambda\right)  $ only satisfies this linear
equation if $\lambda$ is an eigenvalue of $A$ such that $x_{j}=x\left(
\lambda_{j}\right)  $. We have made the dependence on the parameter $\lambda$
explicit and can interpret $\lambda$ as a frequency that ranges continuously
over all real numbers. Invoking the Dirac delta-function $\delta\left(
t\right)  $, we can write%
\begin{align*}
\sum_{j=1}^{N}\left(  x_{j}\right)  _{m}\left(  x_{j}\right)  _{k}  &
=\sum_{\lambda\in\left\{  \lambda_{1},\lambda_{2},\ldots,\lambda_{N}\right\}
}^{N}\left(  x\left(  \lambda\right)  \right)  _{m}\left(  x\left(
\lambda\right)  \right)  _{k}\\
&  =\sum_{j=1}^{N}\int_{-\infty}^{\infty}\delta\left(  \lambda-\lambda
_{j}\right)  \left(  x\left(  \lambda\right)  \right)  _{m}\left(  x\left(
\lambda\right)  \right)  _{k}d\lambda
\end{align*}
Using the non-negative weight function%
\[
w\left(  \lambda\right)  =\sum_{j=1}^{N}\delta\left(  \lambda-\lambda
_{j}\right)  =\delta\left(  \det\left(  A-\lambda I\right)  \right)
\left\vert \left.  \frac{d\det\left(  A-xI\right)  }{dx}\right\vert
_{x=\lambda}\right\vert
\]
shows that%
\begin{equation}
\sum_{j=1}^{N}\left(  x_{j}\right)  _{m}\left(  x_{j}\right)  _{k}%
=\int_{-\infty}^{\infty}w\left(  \lambda\right)  \left(  x\left(
\lambda\right)  \right)  _{m}\left(  x\left(  \lambda\right)  \right)
_{k}d\lambda=\delta_{mk} \label{orthogonality_polynomials}%
\end{equation}
The right-hand side in (\ref{orthogonality_polynomials}) is the continuous
variant of (\ref{second_orthogonality_equations}) that expresses orthogonality
between functions with respect to the weight function~$w$ (see e.g. \cite[p.
313]{PVM_graphspectra}). Specifically\footnote{The eigendecomposition of a
general tri-diagonal stochastic matrix in \cite[Appendix]{PVM_decay_SIS2014}
exemplifies how orthogonal polynomials as a function of $\lambda$ enter.}, the
orthogonality property (\ref{orthogonality_polynomials}) shows that the set
$\left\{  \left(  x\left(  \lambda\right)  \right)  _{m}\right\}  _{1\leq
m\leq N}$ is a set of $N$ orthogonal polynomials in $\lambda$.

\subsection{Calculus for eigenvectors}

Another advantage of the parametrized eigenvalue equation
(\ref{eigenvalue_equation_lambda}) is that calculus can be applied. Invoking
Leibniz' rule, the $n$-th derivative of both sides of $Ax\left(
\lambda\right)  =\lambda\ x\left(  \lambda\right)  $ with respect to $\lambda$
is%
\[
A\frac{d^{n}x\left(  \lambda\right)  }{d\lambda^{n}}=\sum_{k=0}^{n}\binom
{n}{k}\frac{d^{k}}{d\lambda^{k}}\left(  \lambda\right)  \frac{d^{n-k}%
}{d\lambda^{n-k}}x\left(  \lambda\right)  =\lambda\frac{d^{n}x\left(
\lambda\right)  }{d\lambda^{n}}+n\frac{d^{n-1}x\left(  \lambda\right)
}{d\lambda^{n-1}}%
\]
so that, for $n\geq1$,%
\begin{equation}
\left(  A-\lambda I\right)  \frac{d^{n}x\left(  \lambda\right)  }{d\lambda
^{n}}=n\frac{d^{n-1}x\left(  \lambda\right)  }{d\lambda^{n-1}}
\label{Leibniz_derivative}%
\end{equation}
Explicitly, denoting $x^{\left(  n\right)  }\left(  \lambda\right)
=\frac{d^{n}x\left(  \lambda\right)  }{d\lambda^{n}}$, we obtain the sequence%
\begin{align*}
\left(  A-\lambda I\right)  x\left(  \lambda\right)   &  =0\\
\left(  A-\lambda I\right)  x^{\left(  1\right)  }\left(  \lambda\right)   &
=x\left(  \lambda\right) \\
\left(  A-\lambda I\right)  x^{\left(  2\right)  }\left(  \lambda\right)   &
=2x^{\left(  1\right)  }\left(  \lambda\right) \\
&  \vdots\\
\left(  A-\lambda I\right)  x^{\left(  n\right)  }\left(  \lambda\right)   &
=nx^{\left(  n-1\right)  }\left(  \lambda\right)
\end{align*}
from which we deduce that%
\begin{equation}
\left(  A-\lambda I\right)  ^{n+1}x^{\left(  n\right)  }\left(  \lambda
\right)  =0 \label{governing_equation_n_derivative_x(lambda)}%
\end{equation}
but%
\begin{equation}
\left(  A-\lambda I\right)  ^{n}x^{\left(  n\right)  }\left(  \lambda\right)
=n!x\left(  \lambda\right)  \label{x_j_in_terms_x_higher_grade}%
\end{equation}
If $\lambda$ is not an eigenvalue so that $A-\lambda I$ is of rank $N$ and
invertible, then the above shows that $x\left(  \lambda\right)  =0$ (as well
as all higher order derivatives). If $\lambda$ is an eigenvalue, the vector
$x^{\left(  n\right)  }\left(  \lambda\right)  $ can be different from the
zero vector and orthogonal to all the row vectors of $\left(  A-\lambda
I\right)  ^{n+1}$.

\begin{theorem}
\label{theorem_principal_vector_grade_n_independent}The set of vectors
$\left\{  x\left(  \lambda\right)  ,x^{\left(  1\right)  }\left(
\lambda\right)  ,x^{\left(  2\right)  }\left(  \lambda\right)  ,\ldots
,x^{\left(  n\right)  }\left(  \lambda\right)  \right\}  $ is linearly independent.
\end{theorem}

\textbf{Proof:} Assume, on the contrary, that these vectors are dependent,
then%
\[
b_{0}x\left(  \lambda\right)  +b_{1}x^{\left(  1\right)  }\left(
\lambda\right)  +b_{2}x^{\left(  2\right)  }\left(  \lambda\right)
+\ldots+b_{n}x^{\left(  n\right)  }\left(  \lambda\right)  =\sum_{j=0}%
^{n}b_{j}x^{\left(  j\right)  }\left(  \lambda\right)  =0
\]
and not all $b_{j}$ are zero. Left-multiplying both sides with $\left(
A-\lambda I\right)  ^{n}$ and taking into account that $\left(  A-\lambda
I\right)  ^{m+j}x^{\left(  m\right)  }\left(  \lambda\right)  =0$ for any
$j\geq1$ leads to%
\[
b_{n}\left(  A-\lambda I\right)  ^{n}x^{\left(  n\right)  }\left(
\lambda\right)  =0
\]
and (\ref{x_j_in_terms_x_higher_grade}) indicates that $b_{n}$ must be zero.
Next, we repeat the argument and left-multiply both sides with $\left(
A-\lambda I\right)  ^{n-1}$, which leads us to conclude that $b_{n-1}=0$.
Continuing in this way shows that each coefficient $b_{j}=0$ for $0\leq j\leq
n$, which proves the Theorem
\ref{theorem_principal_vector_grade_n_independent}.\hfill$\square\medskip$

Consider for $1\leq k\leq n$ the vectors%
\[
y_{k}=\left(  A-\lambda I\right)  ^{k}x^{\left(  n\right)  }\left(
\lambda\right)
\]
Relation (\ref{x_j_in_terms_x_higher_grade}) shows that $y_{n}=n!x\left(
\lambda\right)  $, while applying (\ref{Leibniz_derivative}) iteratively
$m$-times yields%
\[
y_{k}=\frac{n!}{\left(  n-m\right)  !}\left(  A-\lambda I\right)
^{k-m}x^{\left(  n-m\right)  }\left(  \lambda\right)
\]
from which we find%
\[
y_{m}=\frac{n!}{\left(  n-m\right)  !}x^{\left(  n-m\right)  }\left(
\lambda\right)  =\left(  A-\lambda I\right)  ^{m}x^{\left(  n\right)  }\left(
\lambda\right)
\]
Hence, any vector $y_{m}$ is generated by the vector $x^{\left(  n\right)
}\left(  \lambda\right)  $ and Theorem
\ref{theorem_principal_vector_grade_n_independent} states that the set
$\left\{  y_{1},y_{2},\ldots,y_{n}\right\}  $ is linearly independent and thus
spans the $n$-dimensional space. In the classical eigenvalue theory \cite[p.
43]{Wilkinson}, the vector $z$ satisfying $\left(  A-\lambda_{k}I\right)
^{n+1}z=0$ is called a \emph{principal vector of grade }$n+1$ corresponding to
eigenvalue $\lambda_{k}$. Theorem
\ref{theorem_principal_vector_grade_n_independent} and
(\ref{governing_equation_n_derivative_x(lambda)}) show that $z=\beta
x^{\left(  n\right)  }\left(  \lambda\right)  $, for any non-zero number
$\beta$.

Left-multiplying (\ref{x_j_in_terms_x_higher_grade}) by $x^{T}\left(
\xi\right)  $ yields%
\[
n!x^{T}\left(  \xi\right)  x\left(  \lambda\right)  =x^{T}\left(  \xi\right)
\left(  A-\lambda I\right)  ^{n}x^{\left(  n\right)  }\left(  \lambda\right)
\]
If $A$ is a symmetric matrix and $\xi$ is an eigenvalue of $A$, then
$x^{T}\left(  \xi\right)  \left(  A-\lambda I\right)  ^{n}=x^{T}\left(
\xi\right)  \left(  \xi-\lambda\right)  ^{n}$, so that%
\[
\delta_{\xi\lambda}=x^{T}\left(  \xi\right)  x\left(  \lambda\right)
=\frac{\left(  \xi-\lambda\right)  ^{n}}{n!}x^{T}\left(  \xi\right)
x^{\left(  n\right)  }\left(  \lambda\right)
\]
Hence, if the eigenvalue $\xi$ is different from the eigenvalue $\lambda$, we
find that $x^{T}\left(  \xi\right)  x^{\left(  n\right)  }\left(
\lambda\right)  =0$ for all $n\geq0$. However, when $\xi=\lambda$, an
inconsistency appears when $n>0$, which implies that a \emph{principal vector
of grade }$n+1$ vector $x^{\left(  n\right)  }\left(  \lambda\right)  $ with
$n>0$ does not exist for symmetric matrices. Another argument is that, for
symmetric matrices, the set of eigenvectors $\left\{  x_{m}\right\}  _{1\leq
m\leq N}$ spans the entire space so that $x^{\left(  n\right)  }\left(
\lambda\right)  =0$ for $n\geq1$, because a non-zero vector cannot be
orthogonal to all eigenvectors. Hence, a \emph{principal vector }$x^{\left(
n\right)  }\left(  \lambda\right)  $ \emph{of grade }$n+1$ with $n>0$ can only
exist for asymmetric matrices and may be helpful to construct an orthogonal
set of vectors when degeneracy occurs (as in Jordan forms).

\subsection{Function of a symmetric matrix}

From the general relation for diagonalizable matrices (see e.g. \cite[p.
526]{Meyer_matrix}),
\begin{equation}
f\left(  A\right)  =\sum_{k=1}^{N}f\left(  \lambda_{k}\right)  x_{k}x_{k}%
^{T}\label{eigenvalue_decompostion_function_f}%
\end{equation}
valid for a function $f$ defined on the eigenvalues $\left\{  \lambda
_{k}\right\}  _{1\leq k\leq N}$ of the $N\times N$ matrix $A$, the element for
node $j$ equals%
\begin{equation}
\left(  f\left(  A\right)  \right)  _{jj}=\sum_{k=1}^{N}f\left(  \lambda
_{k}\right)  \left(  x_{k}\right)  _{j}^{2}%
\label{eigenvalue_decompostion_function_f_node_jj}%
\end{equation}
Explicitly, we have%
\[
\left[
\begin{array}
[c]{c}%
\left(  f\left(  A\right)  \right)  _{11}\\
\left(  f\left(  A\right)  \right)  _{22}\\
\left(  f\left(  A\right)  \right)  _{33}\\
\vdots\\
\left(  f\left(  A\right)  \right)  _{NN}%
\end{array}
\right]  =\left[
\begin{array}
[c]{ccccc}%
\left(  x_{1}\right)  _{1}^{2} & \left(  x_{2}\right)  _{1}^{2} & \left(
x_{3}\right)  _{1}^{2} & \cdots & \left(  x_{N}\right)  _{1}^{2}\\
\left(  x_{1}\right)  _{2}^{2} & \left(  x_{2}\right)  _{2}^{2} & \left(
x_{3}\right)  _{2}^{2} & \cdots & \left(  x_{N}\right)  _{2}^{2}\\
\left(  x_{1}\right)  _{3}^{2} & \left(  x_{2}\right)  _{3}^{2} & \left(
x_{3}\right)  _{3}^{2} & \cdots & \left(  x_{N}\right)  _{3}^{2}\\
\vdots & \vdots & \vdots & \ddots & \vdots\\
\left(  x_{1}\right)  _{N}^{2} & \left(  x_{2}\right)  _{N}^{2} & \left(
x_{3}\right)  _{N}^{2} & \cdots & \left(  x_{N}\right)  _{N}^{2}%
\end{array}
\right]  \left[
\begin{array}
[c]{c}%
f\left(  \lambda_{1}\right)  \\
f\left(  \lambda_{2}\right)  \\
f\left(  \lambda_{3}\right)  \\
\vdots\\
f\left(  \lambda_{N}\right)
\end{array}
\right]
\]
which we write in matrix form as $\psi=C\chi$, with%
\[
\psi=\left[
\begin{array}
[c]{c}%
\left(  f\left(  A\right)  \right)  _{11}\\
\left(  f\left(  A\right)  \right)  _{22}\\
\left(  f\left(  A\right)  \right)  _{33}\\
\vdots\\
\left(  f\left(  A\right)  \right)  _{NN}%
\end{array}
\right]  \text{ and }\chi=\left[
\begin{array}
[c]{c}%
f\left(  \lambda_{1}\right)  \\
f\left(  \lambda_{2}\right)  \\
f\left(  \lambda_{3}\right)  \\
\vdots\\
f\left(  \lambda_{N}\right)
\end{array}
\right]
\]
and the matrix $C=X\circ X$, where $\circ$ denotes the Hadamard
product\footnote{The Hadamard product (entrywise product) of two matrix is
$\left(  A\circ B\right)  _{ij}=A_{ij}B_{ij}$. If $A$ and $B$ are both
diagonal matrices, then $A.B=A\circ B$.},%
\begin{equation}
C=\left[
\begin{array}
[c]{ccccc}%
\left(  x_{1}\right)  _{1}^{2} & \left(  x_{2}\right)  _{1}^{2} & \left(
x_{3}\right)  _{1}^{2} & \cdots & \left(  x_{N}\right)  _{1}^{2}\\
\left(  x_{1}\right)  _{2}^{2} & \left(  x_{2}\right)  _{2}^{2} & \left(
x_{3}\right)  _{2}^{2} & \cdots & \left(  x_{N}\right)  _{2}^{2}\\
\left(  x_{1}\right)  _{3}^{2} & \left(  x_{2}\right)  _{3}^{2} & \left(
x_{3}\right)  _{3}^{2} & \cdots & \left(  x_{N}\right)  _{3}^{2}\\
\vdots & \vdots & \vdots & \ddots & \vdots\\
\left(  x_{1}\right)  _{N}^{2} & \left(  x_{2}\right)  _{N}^{2} & \left(
x_{3}\right)  _{N}^{2} & \cdots & \left(  x_{N}\right)  _{N}^{2}%
\end{array}
\right]  \label{def_Y_matrix_squared_eigenvector_components}%
\end{equation}
Since $Cu=u$ and $C^{T}u=u$, by \textquotedblleft double
orthogonality\textquotedblright\ of (\ref{first_orthogonality_equations}) and
(\ref{second_orthogonality_equations}), and since each element $0\leq\left(
x_{k}\right)  _{j}^{2}\leq$ $1$, the matrix $Y$ with squared eigenvector
components of a diagonalizable matrix $A$ is doubly\footnote{\emph{Sinkhorn's
theorem} (1964) states that any matrix with strictly positive entries can be
made doubly stochastic by pre- and post-multiplication by diagonal matrices.}
stochastic \cite{PVM_graphspectra} with largest eigenvalue equal to 1.

Let us denote the vector $\lambda^{k}=\left(  \lambda_{1}^{k},\lambda_{2}%
^{k},\ldots,\lambda_{N}^{k}\right)  $ so that, for $f\left(  z\right)  =z^{n}%
$, we have%
\begin{equation}
\text{diag}\left(  \left(  A^{k}\right)  _{jj}\right)  u=C\lambda
^{k}\label{general_diag_A^k_in_Y}%
\end{equation}
From (\ref{general_diag_A^k_in_Y}) and $u^{T}C=u^{T}$, we find the well-known
trace relation, namely that $u^{T}$diag$\left(  \left(  A^{k}\right)
_{jj}\right)  u=$ trace$\left(  A^{k}\right)  =u^{T}\lambda^{k}=\sum_{j=1}%
^{N}\lambda_{j}^{k}$. If the inverse $C^{-1}$ of $C$ exists, then it holds,
for any integer $k$, that%
\[
\lambda^{k}=C^{-1}\text{diag}\left(  \left(  A^{k}\right)  _{jj}\right)  u
\]
or, the eigenvalue $\lambda_{j}$ of $A$ (to any integer power $k$) can be
written as a linear combination of the diagonal elements of $A^{k}$,%
\[
\lambda_{j}^{k}=\sum_{i=1}^{N}\left(  C^{-1}\right)  _{ji}\left(
A^{k}\right)  _{ii}%
\]

We can proceed on step further by applying the above to a set $f_{1}%
,f_{2},\ldots,f_{N}$ of $N$ functions and obtain the matrix equation%
\[
F=CG
\]
where the $N\times N$ matrix $F$ is%
\[
F=\left[
\begin{array}
[c]{ccccc}%
\left(  f_{1}\left(  A\right)  \right)  _{11} & \left(  f_{2}\left(  A\right)
\right)  _{11} & \left(  f_{3}\left(  A\right)  \right)  _{11} & \cdots &
\left(  f_{N}\left(  A\right)  \right)  _{11}\\
\left(  f_{1}\left(  A\right)  \right)  _{22} & \left(  f_{2}\left(  A\right)
\right)  _{22} & \left(  f_{3}\left(  A\right)  \right)  _{22} & \cdots &
\left(  f_{N}\left(  A\right)  \right)  _{22}\\
\left(  f_{1}\left(  A\right)  \right)  _{33} & \left(  f_{2}\left(  A\right)
\right)  _{33} & \left(  f_{3}\left(  A\right)  \right)  _{33} & \cdots &
\left(  f_{N}\left(  A\right)  \right)  _{33}\\
\vdots & \vdots & \vdots & \ddots & \vdots\\
\left(  f_{1}\left(  A\right)  \right)  _{NN} & \left(  f_{2}\left(  A\right)
\right)  _{NN} & \left(  f_{3}\left(  A\right)  \right)  _{NN} & \cdots &
\left(  f_{N}\left(  A\right)  \right)  _{NN}%
\end{array}
\right]
\]
and the $N\times N$ matrix $G$ is%
\[
G=\left[
\begin{array}
[c]{ccccc}%
f_{1}\left(  \lambda_{1}\right)   & f_{2}\left(  \lambda_{1}\right)   &
f_{3}\left(  \lambda_{1}\right)   & \cdots & f_{N}\left(  \lambda_{1}\right)
\\
f_{1}\left(  \lambda_{2}\right)   & f_{2}\left(  \lambda_{2}\right)   &
f_{3}\left(  \lambda_{2}\right)   & \cdots & f_{N}\left(  \lambda_{2}\right)
\\
f_{1}\left(  \lambda_{3}\right)   & f_{2}\left(  \lambda_{3}\right)   &
f_{3}\left(  \lambda_{3}\right)   & \cdots & f_{N}\left(  \lambda_{3}\right)
\\
\vdots & \vdots & \vdots & \ddots & \vdots\\
f_{1}\left(  \lambda_{N}\right)   & f_{2}\left(  \lambda_{N}\right)   &
f_{3}\left(  \lambda_{N}\right)   & \cdots & f_{N}\left(  \lambda_{N}\right)
\end{array}
\right]
\]
If $G$ is invertible (i.e. $\det G\neq0$), which requires that all eigenvalues
are distinct, then we can construct $C=FG^{-1}$ from which we deduce that
$\det Y=\frac{\det F}{\det G}$. A straightforward choice are the functions
$f_{n}\left(  z\right)  =z^{n-1}$, so that $G$ reduces to a Vandermonde
matrix, in which case, $C=FG^{-1}$ leads to the results in Theorem
\ref{theorem_eigenvectors_in_terms_walks}.

\subsubsection{Application to the adjacency matrix}

\label{sec_centrality_vectors_are_dependent}Applied to the adjacency matrix
$A$, (\ref{eigenvalue_decompostion_function_f_node_jj}) illustrates that the
\emph{squares} of the eigenvector component arise as weights for $f\left(
\lambda_{k}\right)  $ to specify a function of the adjacency matrix $A$ at
node $j$. In particular, for powers $f\left(  z\right)  =z^{n}$, nice formulae
appear: for $n=0$, we find from (\ref{eigenvalue_decompostion_function_f}) the
\emph{second} \cite{PVM_double_orthogonality} orthogonality relation
(\ref{second_orthogonality_equations}); for $n=1$ (since $A_{jj}=0$, from
which trace$\left(  A\right)  =\sum_{k=1}^{N}\lambda_{k}=0$)
\begin{equation}
0=\sum_{k=1}^{N}\lambda_{k}\left(  x_{k}\right)  _{j}^{2}\text{ and
}0=C\lambda\label{zero_node_j_square_eigenvector_components}%
\end{equation}
while for $n=2$ (since $\left(  A^{2}\right)  _{jj}=d_{j}$)%
\begin{equation}
d_{j}=\sum_{k=1}^{N}\lambda_{k}^{2}\left(  x_{k}\right)  _{j}^{2}\text{ and
}d=C\lambda^{2}\label{degree_node_j_square_eigenvector_components}%
\end{equation}
For any adjacency matrix $A$, (\ref{zero_node_j_square_eigenvector_components}%
) shows \cite[p. 229]{PVM_graphspectra} that
\begin{equation}
C\lambda=0\label{zero_eigenvalue_equation_Y}%
\end{equation}
so that det$C=0$ and that the vector $\lambda=\left(  \lambda_{1},\lambda
_{2},\ldots,\lambda_{N}\right)  $ is the eigenvector of $C$ corresponding to
eigenvalue zero. Relation (\ref{zero_eigenvalue_equation_Y}) implies that the
centrality vector $c_{j}=$ $\left(  \left(  x_{1}\right)  _{j}^{2},\left(
x_{2}\right)  _{j}^{2},\ldots,\left(  x_{N}\right)  _{j}^{2}\right)  $ is not
independent from $c_{l}$. In other words, at least one vector $c_{k}$ can be
written as a linear combination of all the other nodal centrality vectors
$\left\{  c_{j}\right\}  _{1\leq j\neq k\leq N}$ and the set $\left\{
c_{j}\right\}  _{1\leq j\leq N}$ is not complete, in that it does not span the
entire $N$-dimensional space.  

The fact that $\det C=\det C^{T}$, implies that the (left)-eigenvector $q$ of
$C^{T}$ belonging to the zero eigenvalue obeys, for each $1\leq k\leq N$,%
\[
0=\sum_{j=1}^{N}q_{j}\left(  x_{k}\right)  _{j}^{2}%
\]
which is the companion of (\ref{zero_node_j_square_eigenvector_components})
over the node labels $j$.

\section{Walk expansion}

\label{sec_walk_based_(xk)_j}The following theorem is a direct consequence of
the analysis in \cite[p. 228]{PVM_graphspectra}:

\begin{theorem}
\label{theorem_eigenvectors_in_terms_walks}If all eigenvalues of $A$ are
different, then
\begin{equation}
\left(  x_{k}\right)  _{i}\left(  x_{k}\right)  _{j}=\frac{1}{{\prod
_{l=1;l\neq k}^{N}}\left(  \lambda_{k}-\lambda_{l}\right)  }\sum_{r=H_{ij}%
}^{N-1}b_{r}\left(  k\right)  \left(  A^{r}\right)  _{ij}
\label{prod(xm)_i(xm)_j_all_hops}%
\end{equation}
where $H_{ij}$ is the hopcount (number of links) of the shortest path between
node $i$ and $j$ and where the coefficients $b_{r}\left(  k\right)  $ obey
\[
\prod_{j=1;j\neq k}^{N}\left(  x-\lambda_{j}\right)  =\sum_{j=0}^{N-1}%
b_{j}\left(  k\right)  x^{j}%
\]
or%
\begin{equation}
b_{r}\left(  k\right)  =\frac{1}{r!}\left.  \frac{d^{r}}{dx^{r}}%
\prod_{j=1;j\neq k}^{N}\left(  x-\lambda_{j}\right)  \right\vert _{x=0}
\label{def_b_r(m)}%
\end{equation}

Writing (\ref{prod(xm)_i(xm)_j_all_hops}) in matrix form yields
\begin{equation}
x_{k}x_{k}^{T}=\frac{1}{{\prod_{l=1;l\neq k}^{N}}\left(  \lambda_{k}%
-\lambda_{l}\right)  }\sum_{r=0}^{N-1}b_{r}\left(  k\right)  A^{r}=\frac
{\prod_{j=1;j\neq k}^{N}\left(  A-\lambda_{j}I\right)  }{{\prod_{l=1;l\neq
k}^{N}}\left(  \lambda_{k}-\lambda_{l}\right)  } \label{matrixE_m_in_terms_A}%
\end{equation}

\end{theorem}

Clearly, if $i=j$, then $H_{jj}=0$ and (\ref{prod(xm)_i(xm)_j_all_hops})
reduces to%
\begin{equation}
\left(  x_{k}\right)  _{j}^{2}=\frac{1}{{\prod_{l=1;l\neq k}^{N}}\left(
\lambda_{k}-\lambda_{l}\right)  }\sum_{r=0}^{N-1}\left(  A^{r}\right)
_{jj}b_{r}\left(  k\right)  \label{(x_m)^2}%
\end{equation}
The definition of the characteristic polynomial of matrix $A$ is $c_{A}\left(
\lambda\right)  =\det\left(  A-\lambda I\right)  =\prod_{j=1}^{N}\left(
\lambda_{j}-\lambda\right)  $, from which $\log c_{A}\left(  \lambda\right)
=\sum_{j=1}^{N}\log\left(  \lambda_{j}-\lambda\right)  $. Differentiation
yields%
\[
c_{A}^{\prime}\left(  \lambda\right)  =-c_{A}\left(  \lambda\right)
\sum_{j=1}^{N}\frac{1}{\lambda_{j}-\lambda}=-\sum_{j=1}^{N}\frac{\prod
_{k=1}^{N}\left(  \lambda_{k}-\lambda\right)  }{\lambda_{j}-\lambda}%
=-\sum_{j=1}^{N}\prod_{k=1;k\neq j}^{N}\left(  \lambda_{k}-\lambda\right)
\]
from which%
\begin{equation}
c_{A}^{\prime}\left(  \lambda_{m}\right)  =-\prod_{k=1;k\neq m}^{N}\left(
\lambda_{k}-\lambda_{m}\right)  =\left(  -1\right)  ^{N}\prod_{k=1;k\neq
m}^{N}\left(  \lambda_{m}-\lambda_{k}\right)
\label{afgeleide_characteristic_pol_in_lambda_m}%
\end{equation}
Thus, we can write (\ref{(x_m)^2}) as%
\begin{equation}
\left(  x_{k}\right)  _{j}^{2}=\frac{\left(  -1\right)  ^{N}}{c_{A}^{\prime
}\left(  \lambda_{k}\right)  }\sum_{r=0}^{N-1}\left(  A^{r}\right)  _{jj}%
b_{r}\left(  k\right)  \label{(x_m)^2_char_pol}%
\end{equation}
Theorem \ref{theorem_eigenvectors_in_terms_walks} expresses the product of two
eigenvector components in terms of the eigenvalues only. In particular,
(\ref{prod(xm)_i(xm)_j_all_hops}) equals the sum of the number $\left(
A^{r}\right)  _{ij}$ of walks, weighted by a function $b_{r}\left(  m\right)
$ of eigenvalues, over all $r$ hops paths between node $i$ and node $j$. The
longest path in a graph contains $N-1$ hops and $\left(  A^{r}\right)  _{ij}$
equals the number of \emph{shortest paths} with $r$ hops from node $i$ to node
$j$, provided $\left(  A^{m}\right)  _{ij}=0$ for all integers $m<r$. Invoking
the normalization $x_{k}^{T}x_{k}=\sum_{j=1}^{N}\left(  x_{k}\right)  _{j}%
^{2}=1$ and $W_{r}=\sum_{j=1}^{N}\left(  A^{r}\right)  _{jj}$, the total
number of closed walks of length $r$ (with $r$ hops), we obtain from
(\ref{(x_m)^2_char_pol}) that%
\[
c_{A}^{\prime}\left(  \lambda_{k}\right)  =\left(  -1\right)  ^{N}\sum
_{r=0}^{N-1}W_{r}b_{r}\left(  k\right)
\]

\section{Spacings of vector components}

\label{sec_spacings}We derive bounds for the minimum and maximum spacing in a
vector $b$, whose components are ordered.

\begin{theorem}
\label{theo_spacing_sequence}Let the set of real numbers\ $\left\{
b_{j}\right\}  _{1\leq j\leq n}$ be ordered as $b_{1}\geq b_{2}\geq\cdots\geq
b_{n}$, then for non-negative real numbers $\left\{  a_{j}\right\}  _{1\leq
j\leq n}$ with at least one $a_{j}>0$, we have that%
\begin{equation}
\min_{1\leq j\leq n}\left(  b_{j}-b_{j+1}\right)  \leq\frac{\sum_{k=1}%
^{n-1}a_{k}\,\left(  b_{k}-b_{n}\right)  }{n\sum_{l=1}^{n}a_{l}-\sum_{l=1}%
^{n}la_{l}}\leq\max_{1\leq j\leq n}\left(  b_{j}-b_{j+1}\right)
\label{general_spacing_inequality}%
\end{equation}
with equality only if all spacings are the same.\qquad
\end{theorem}

\textbf{Proof:} We rewrite Abel's summation formula \cite[p. 56]%
{PVM_graphspectra} as%
\[
\sum_{k=1}^{n-1}a_{k}\,\left(  b_{k}-b_{n}\right)  =\sum_{k=1}^{n-1}\left(
\sum_{l=1}^{k}a_{l}\right)  \,(b_{k}-b_{k+1})
\]
Since $b_{k}-b_{k+1}\geq0$ and $\sum_{l=1}^{k}a_{l}\geq0$, we bound the
right-hand side as%
\[
\min_{1\leq j\leq N}\left(  b_{j}-b_{j+1}\right)  \sum_{k=1}^{n-1}\left(
\sum_{l=1}^{k}a_{l}\right)  \leq\sum_{k=1}^{n-1}\left(  \sum_{l=1}^{k}%
a_{l}\right)  \,(b_{k}-b_{k+1})\leq\max_{1\leq j\leq N}\left(  b_{j}%
-b_{j+1}\right)  \sum_{k=1}^{n-1}\left(  \sum_{l=1}^{k}a_{l}\right)
\]
Further,%
\[
0\leq\sum_{k=1}^{n-1}\left(  \sum_{l=1}^{k}a_{l}\right)  =\sum_{l=1}^{n-1}%
\sum_{k=l}^{n-1}a_{l}=\sum_{l=1}^{n-1}\left(  n-l\right)  a_{l}=n\sum
_{l=1}^{n-1}a_{l}-\sum_{l=1}^{n-1}la_{l}=n\sum_{l=1}^{n}a_{l}-\sum_{l=1}%
^{n}la_{l}%
\]
Combining all leads to (\ref{general_spacing_inequality}).\hfill
$\square\medskip$

Theorem \ref{theo_spacing_sequence} illustrates that spacing of the ordered
series $\left\{  b_{j}\right\}  _{1\leq j\leq n}$ are the same as that of
$\left\{  c_{j}\right\}  _{1\leq j\leq n}$, where $c_{j}=b_{j}-b_{n}$ and
$c_{n}=0$. Indeed, $b_{k}-b_{k+1}=\left(  b_{k}-b_{n}\right)  -\left(
b_{k+1}-b_{n}\right)  =c_{k}-c_{k+1}$. If we denote the $\left(  n-1\right)
\times1$ vectors $a=\left(  a_{1},a_{2},\ldots,a_{n-1}\right)  $, $b=\left(
b_{1},b_{2},\ldots,b_{n-1}\right)  $ and $t=\left(  n-1,n-2,n-3,\ldots
,1\right)  $, so that $c=b-\left(  b_{n}\right)  u$, then the fraction in
Theorem \ref{theo_spacing_sequence} can be written as%
\begin{equation}
f=\frac{\sum_{k=1}^{n-1}a_{k}\,\left(  b_{k}-b_{n}\right)  }{\sum_{l=1}%
^{n-1}\left(  n-l\right)  a_{l}}=\frac{a^{T}c}{a^{T}t}=\frac{\left\Vert
c\right\Vert _{2}}{\left\Vert t\right\Vert _{2}}\frac{\left(  \frac
{a}{\left\Vert a\right\Vert _{2}}\right)  ^{T}\frac{c}{\left\Vert c\right\Vert
_{2}}}{\left(  \frac{a}{\left\Vert a\right\Vert _{2}}\right)  ^{T}\frac
{t}{\left\Vert t\right\Vert _{2}}} \label{def_f}%
\end{equation}
where in the last equality, $f$ is written terms of normalized vectors, where
$\left\Vert t\right\Vert _{2}^{2}=\sum_{l=1}^{n-1}\left(  n-l\right)
^{2}=\sum_{k=1}^{n-1}k^{2}=\frac{\left(  n-1\right)  n\left(  2n-1\right)
}{6}$. If $a=c$, then%
\[
f=\frac{\left\Vert c\right\Vert _{2}^{2}}{c^{T}t}\geq\frac{\left\Vert
c\right\Vert _{2}}{\left\Vert t\right\Vert _{2}}%
\]
If $a=e_{n-1}$, then%
\[
f=b_{n-1}-b_{n}%
\]
while if $a=e_{1}$, then $f=\frac{b_{1}-b_{n}}{n-1}$ and
(\ref{general_spacing_inequality}) reduces to bounds deduced from the
telescoping series $\sum_{j=1}^{n-1}\left(  b_{j}-b_{j+1}\right)  =b_{1}%
-b_{n}$. If $a=e_{1}-e_{2}$, then $f=b_{1}-b_{2}$. Finally, if $a=u$, then
(\ref{def_f}) reduces, with $u^{T}t=\binom{n}{2}$ to
\[
f=\frac{u^{T}c}{\binom{n}{2}}=\frac{2\left(  \frac{1}{n}\sum_{k=1}^{n}%
b_{k}-b_{n}\right)  }{n-1}%
\]
Finding the vector $a$ that either maximizes or minimizes $f$ would be useful.

\section{Eigenvectors of the complementary graph}

\label{sec_eigenvectors_complementary_graph}The adjacency matrix of the
complementary graph $G^{c}$ is $A^{c}=J-I-A$. In general, $A^{c}$ does not
commute with $A$, unless the graph is regular \cite[p. 44]{PVM_graphspectra}.
When symmetric matrices commute, the eigenvectors are the same. Let $z_{k}$ be
the eigenvector of $A^{c}$ belonging to eigenvalue $\theta_{k}$, so that
$A^{c}z_{k}=\theta_{k}z_{k}$. Since both $A$ and $A^{c}$ are symmetric, a
complete set of eigenvectors exists, so that%
\[
z_{k}=\sum_{j=1}^{N}\left(  z_{k}^{T}x_{j}\right)  x_{j}\text{ and }x_{m}%
=\sum_{j=1}^{N}\left(  x_{m}^{T}z_{j}\right)  z_{j}%
\]
Now,%
\begin{align*}
A^{c}x_{m}  &  =Jx_{m}-\left(  \lambda_{m}+1\right)  x_{m}\\
&  =w_{m}u-\left(  \lambda_{m}+1\right)  x_{m}%
\end{align*}
where, $Jx_{m}=u.u^{T}x_{m}=w_{m}u$ and similarly, $Jz_{k}=v_{k}u$, where
$v_{k}=u^{T}z_{k}$, the fundamental weight of the complementary graph $G^{c}$.
Left-multiplying $A^{c}z_{k}$ by $x_{m}^{T}$ and $A^{c}x_{m}$ by $z_{k}^{T}$
yields
\[
x_{m}^{T}A^{c}z_{k}=\theta_{k}x_{m}^{T}z_{k}%
\]
and%
\[
z_{k}^{T}A^{c}x_{m}=w_{m}v_{k}-\left(  \lambda_{m}+1\right)  z_{k}^{T}x_{m}%
\]
Since $z_{m}^{T}A^{c}y_{k}=z_{k}^{T}A^{c}x_{m}$ (because it is a scalar), we
deduce that%
\[
x_{m}^{T}z_{k}=\frac{w_{m}v_{k}}{\theta_{k}+\lambda_{m}+1}%
\]
and%
\[
z_{k}=v_{k}\sum_{j=1}^{N}\frac{w_{j}}{\theta_{k}+\lambda_{j}+1}x_{j}\text{ and
}x_{m}=w_{m}\sum_{j=1}^{N}\frac{v_{j}}{\theta_{j}+\lambda_{m}+1}z_{j}%
\]
The last, when left-multiplied with $u^{T}$ yields, for any $k$,%
\[
\sum_{j=1}^{N}\frac{w_{j}^{2}}{\theta_{k}+\lambda_{j}+1}=1
\]
and similarly\footnote{Using the orthogonality relations $z_{k}^{T}%
z_{m}=\delta_{km}$, these equations are complemented by%
\[
\sum_{j=1}^{N}\frac{w_{j}^{2}}{\left(  \theta_{k}+\lambda_{j}+1\right)
\left(  \theta_{m}+\lambda_{j}+1\right)  }=\frac{\delta_{km}}{v_{k}v_{m}}%
\]
and, similarly,%
\[
\sum_{j=1}^{N}\frac{v_{j}^{2}}{\left(  \theta_{k}+\lambda_{j}+1\right)
\left(  \theta_{m}+\lambda_{j}+1\right)  }=\frac{\delta_{km}}{w_{k}w_{m}}%
\]
}, for any $m$,
\[
\sum_{j=1}^{N}\frac{v_{j}^{2}}{\theta_{j}+\lambda_{m}+1}=1
\]
while, in general, $\sum_{j=1}^{N}w_{j}^{2}=N$. Invoking
(\ref{inequality_breuk_n_termen}) yields, for any $k$,%
\[
\min_{1\leq j\leq N}\left(  \theta_{k}+\lambda_{j}+1\right)  \leq N\leq
\max_{1\leq j\leq N}\left(  \theta_{k}+\lambda_{j}+1\right)
\]
Thus, any eigenvalue $\theta_{k}$ of the complementary adjacency matrix
$A^{c}$ can be bounded in terms of eigenvalues $\lambda_{j}$ of $A$. For
example, for $\theta_{1}>0$,
\[
0\leq N-1-\lambda_{1}\leq\theta_{1}\leq N-1-\lambda_{N}%
\]
where the upper bound is useless. We cannot use $2L=\sum_{j=1}^{N}\lambda
_{j}w_{j}^{2}$ to derive sharper bounds, because all terms must be positive
for the denominator of (\ref{inequality_breuk_n_termen}), but we can use
$N_{2}=d^{T}d=\sum_{j=1}^{N}\lambda_{j}^{2}w_{j}^{2}$. However, we rather
prefer to follow another track by computing $x_{m}^{T}\left(  A^{c}\right)
^{n}z_{k}=\theta_{k}^{n}x_{m}^{T}z_{k}$ for any integer $n\geq1$ as $z_{k}%
^{T}\left(  A^{c}\right)  ^{n}x_{m}$ using%
\[
\left(  A^{c}\right)  ^{n}x_{m}=\left(  J-I-A\right)  ^{n}x_{m}%
\]
After some tedious computations, we find%
\[
\left(  J-I-A\right)  ^{n}x_{m}=x_{m}\left(  -1\right)  ^{n}\left(
\lambda_{m}+1\right)  ^{n}+\frac{w_{m}u}{N}\left\{  \left(  N-1-\lambda
_{m}\right)  ^{n}-\left(  -\lambda_{m}-1\right)  ^{n}\right\}
\]
For example, for $n=1$, we find again the above. Equating $z_{k}^{T}\left(
A^{c}\right)  ^{n}x_{m}=$ $x_{m}^{T}\left(  A^{c}\right)  ^{n}z_{k}$ yields,
for any integer $n\geq1$,%
\[
x_{m}^{T}z_{k}=\frac{w_{m}v_{k}}{N}\frac{\left(  N-1-\lambda_{m}\right)
^{n}-\left(  -1\right)  ^{n}\left(  \lambda_{m}+1\right)  ^{n}}{\theta_{k}%
^{n}-\left(  -1\right)  ^{n}\left(  \lambda_{m}+1\right)  ^{n}}%
\]
Finally, we find the generalized expression for the eigenvector $z_{k}$ of
$G^{c}$ in terms of those of $G$,%
\[
z_{k}=\frac{v_{k}}{N}\sum_{j=1}^{N}\frac{\left(  N-1-\lambda_{j}\right)
^{n}-\left(  -1\right)  ^{n}\left(  \lambda_{j}+1\right)  ^{n}}{\theta_{k}%
^{n}-\left(  -1\right)  ^{n}\left(  \lambda_{j}+1\right)  ^{n}}w_{j}x_{j}%
\]
and, vice versa,%
\[
x_{m}=\frac{w_{m}}{N}\sum_{j=1}^{N}\frac{\left(  N-1-\lambda_{m}\right)
^{n}-\left(  -1\right)  ^{n}\left(  \lambda_{m}+1\right)  ^{n}}{\theta_{j}%
^{n}-\left(  -1\right)  ^{n}\left(  \lambda_{m}+1\right)  ^{n}}v_{j}z_{j}%
\]
After multiplying both sides with $u^{T}$, it follows that%
\[
N=\sum_{j=1}^{N}\frac{\left(  N-1-\lambda_{j}\right)  ^{n}-\left(  -1\right)
^{n}\left(  \lambda_{j}+1\right)  ^{n}}{\theta_{k}^{n}-\left(  -1\right)
^{n}\left(  \lambda_{j}+1\right)  ^{n}}w_{j}^{2}%
\]
Similarly as above, invoking (\ref{inequality_breuk_n_termen}) yields, for any
$1\leq k\leq N$ and $n\geq1$,%
\[
\min_{1\leq j\leq N}\frac{1-\left(  1-\frac{N}{1+\lambda_{j}}\right)  ^{n}%
}{1-\left(  -\frac{\theta_{k}}{1+\lambda_{j}}\right)  ^{n}}\leq1\leq
\max_{1\leq j\leq N}\frac{1-\left(  1-\frac{N}{1+\lambda_{j}}\right)  ^{n}%
}{1-\left(  -\frac{\theta_{k}}{1+\lambda_{j}}\right)  ^{n}}%
\]
This inequality can be used to derive bounds for any eigenvalue $\theta_{k}$
of $A^{c}$ in terms of eigenvalues $\lambda_{j}$ of $A$ by optimizing $n$. The
presented approach complements the determinant theory of $\det\left(
A^{c}-\lambda I\right)  $ in \cite[p. 42-43]{PVM_graphspectra}.

\section{Additions to Theorem \ref{theorem_eigenvector_component_determinants}%
}

\label{sec_addition_Theorem_square_xk_determinant}

\subsection{Introducing the resolvent}

Another way to rewrite the determinant in (\ref{def_scaling_alfa_m}) is%
\[
\det\left(  A-\lambda_{k}I\right)  _{\operatorname{row}N=b}=\det\left[
\begin{array}
[c]{cc}%
\left(  A_{G\backslash\left\{  N\right\}  }-\lambda_{k}I\right)  & \left(
a_{N}\right)  _{\backslash N}\\
b_{\backslash N}^{T} & b_{N}%
\end{array}
\right]
\]
where the $\left(  N-1\right)  \times1$ vector $w_{\backslash m}=\left(
w_{1},\ldots,w_{m-1},w_{m+1},\ldots,w_{N}\right)  $ is obtained from the
$N\times1$ vector $w$ after removing the $m$-th component. Invoking Schur's
block determinant relation \cite[p. 255]{PVM_graphspectra} yields\footnote{We
remark that, in case $b=u$, then
\[
\det\left(  A_{G_{\text{cone}\left(  N\right)  }}-\lambda I\right)
=\det\left[
\begin{array}
[c]{cc}%
\left(  A_{G\backslash\left\{  N\right\}  }-\lambda I\right)  & u\\
u^{T} & -\lambda
\end{array}
\right]
\]
where $G_{\text{cone}\left(  j\right)  }$ is the \textquotedblleft cone at
node $j$\textquotedblright\ of the original graph $G$, which is the graph
where only node $j$ has now links to all other nodes in $G$. In other words,
the node $j$ is the cone of the graph $G\backslash\left\{  j\right\}  $. Thus,
even if $a_{N}=u$, $\det\left(  A-\lambda I\right)  _{\operatorname{row}N=u}$
is not equal to $\det\left(  A_{G_{\text{cone}\left(  N\right)  }}-\lambda
I\right)  $, unless $\lambda=-1$.}%
\[
\det\left[
\begin{array}
[c]{cc}%
\left(  A_{G\backslash\left\{  N\right\}  }-\lambda_{k}I\right)  & \left(
a_{N}\right)  _{\backslash N}\\
b_{\backslash N}^{T} & b_{N}%
\end{array}
\right]  =\det\left(  A_{G\backslash\left\{  N\right\}  }-\lambda_{k}I\right)
\left(  b_{N}-b_{\backslash N}^{T}\left(  A_{G\backslash\left\{  N\right\}
}-\lambda_{k}I\right)  ^{-1}\left(  a_{N}\right)  _{\backslash N}\right)
\]
Instead of row $N$, we can delete row $m$ so that%
\begin{equation}
\det\left(  A-\lambda_{k}I\right)  _{\operatorname{row}m=b}=\det\left(
A_{G\backslash\left\{  m\right\}  }-\lambda_{k}I\right)  \left(
b_{m}-b_{\backslash m}^{T}\left(  A_{G\backslash\left\{  m\right\}  }%
-\lambda_{k}I\right)  ^{-1}\left(  a_{m}\right)  _{\backslash m}\right)
\label{cofactor_expansion_row_m}%
\end{equation}
where $\left(  a_{m}\right)  _{\backslash m}=\left(  a_{1m},\ldots
,a_{m-1;m},a_{m+1,m},\ldots,a_{Nm}\right)  $. Using
(\ref{cofactor_expansion_row_m}) in (\ref{def_scaling_alfa_m}) transforms
(\ref{eigenvector_component_xk_j}) to%
\begin{equation}
\left(  x_{k}\right)  _{j}=\frac{\beta_{k}}{b_{j}-b_{\backslash j}^{T}\left(
A_{G\backslash\left\{  j\right\}  }-\lambda_{k}I\right)  ^{-1}\left(
a_{j}\right)  _{\backslash j}} \label{x_k_j_resolvent}%
\end{equation}
which illustrates the seemingly dependence of $\left(  x_{k}\right)  _{j}$ on
the arbitrary vector $b$.

If $b=e_{m}$, the basic vector with all zero components, except that the
$m$-th component is 1, then (\ref{x_k_j_resolvent}) reduces, for $m\neq j$, to%
\[
\left(  x_{k}\right)  _{j}=-\frac{\left(  x_{k}\right)  _{m}}{\left(  \left(
A_{G\backslash\left\{  j\right\}  }-\lambda_{k}I\right)  ^{-1}\left(
a_{j}\right)  _{\backslash j}\right)  _{m}}%
\]
else, for $m=j$, we find an identity. Interchanging $m$ and $j$, the ratio
$\frac{\left(  x_{k}\right)  _{j}}{\left(  x_{k}\right)  _{m}}$, expressed in
two ways, leads to%
\[
\left(  \left(  A_{G\backslash\left\{  m\right\}  }-\lambda_{k}I\right)
^{-1}\left(  a_{j}\right)  _{\backslash m}\right)  _{j}=\frac{1}{\left(
\left(  A_{G\backslash\left\{  j\right\}  }-\lambda_{k}I\right)  ^{-1}\left(
a_{j}\right)  _{\backslash j}\right)  _{m}}%
\]
When the vector $b$ equals a row vector in $A$, it can be shown (see e.g.
\cite{Tao_Vu2011}) that
\[
\left(  x_{k}\right)  _{j}^{2}=\frac{1}{1+\left(  a_{j}\right)  _{\backslash
j}^{T}\left(  A_{G\backslash\left\{  j\right\}  }-\lambda_{k}I\right)
^{-2}\left(  a_{j}\right)  _{\backslash j}}%
\]

\subsection{Expressions for $\beta_{k}$}

After multiplying (\ref{x_k_j_detA/j}) by $b_{j}$ and summing over all $j$ and
using $\beta_{k}=b^{T}x_{k}=\sum_{j=1}^{N}b_{j}\left(  x_{k}\right)  _{j}$, we
obtain a normalization formula for all $\lambda_{k}$,%
\begin{equation}
\sum_{j=1}^{N}\frac{b_{j}\det\left(  A_{G\backslash\left\{  j\right\}
}-\lambda_{k}I\right)  }{\det\left(  A-\lambda_{k}I\right)
_{\operatorname{row}j=b}}=1 \label{Sum=1}%
\end{equation}
Similarly from (\ref{x_k_j_resolvent}), we obtain%
\begin{equation}
\sum_{j=1}^{N}\frac{b_{j}}{b_{j}-b_{\backslash j}^{T}\left(  A_{G\backslash
\left\{  j\right\}  }-\lambda_{k}I\right)  ^{-1}\left(  a_{j}\right)
_{\backslash j}}=1 \label{Sum=1_resolvent}%
\end{equation}
and\footnote{Substituting (\ref{eigenvector_xk_component_j}) into the
eigenvalue equation \ $\sum_{j=1}^{N}a_{ij}\left(  x_{k}\right)  _{j}%
=\lambda_{k}\left(  x_{k}\right)  _{i}$ gives%
\[
\det\left(  A-\lambda_{k}I\right)  _{\operatorname{row}i=b}=\frac{1}%
{\lambda_{k}}\sum_{j=1}^{N}a_{ij}\det\left(  A-\lambda_{k}I\right)
_{\operatorname{row}j=b}%
\]
} from (\ref{eigenvector_xk_component_j}),%
\begin{equation}
\beta_{k}^{2}=-\frac{1}{c_{A}^{\prime}\left(  \lambda_{k}\right)  }\sum
_{j=1}^{N}b_{j}\det\left(  A-\lambda_{k}I\right)  _{\operatorname{row}j=b}
\label{beta_k^2}%
\end{equation}
After combining (\ref{alfa_m_squared_inverse}) and
(\ref{derivative_characteristic_polynomial_A}) with the definition
(\ref{def_scaling_alfa_m}) of $\alpha_{m}\left(  k\right)  $, we
obtain\footnote{If we choose $b=e_{m}$, then $\beta_{k}=\left(  x_{k}\right)
_{m}$ and $\det\left(  A-\lambda_{k}I\right)  _{\operatorname{row}j=e_{m}%
}=\left(  -1\right)  ^{j+m}\det\left(  A_{G\backslash\operatorname{row}%
j\backslash\operatorname{col}m}-\lambda I\right)  $. Invoking Jacobi's formula
(\ref{Jacobi_generalized_cofactor_theorem_k=2}) in (\ref{beta_k^2_node_m})
leads to (\ref{squared_eigenvector_component_xk_j_charc_pol}).}, for any node
$m$,%
\begin{equation}
\beta_{k}^{2}=-\frac{\det\left(  A-\lambda_{k}I\right)  _{\operatorname{row}%
m=b}^{2}}{c_{A}^{\prime}\left(  \lambda_{k}\right)  \det\left(  A_{G\backslash
\left\{  m\right\}  }-\lambda_{k}I\right)  } \label{beta_k^2_node_m}%
\end{equation}
Summing (\ref{beta_k^2_node_m}) over all $m$, or similarly introducing
(\ref{eigenvector_xk_component_j}) in the first orthogonality relation
$x_{k}^{T}x_{k}=1$, yields, with (\ref{derivative_characteristic_polynomial_A}%
),%
\begin{equation}
\beta_{k}^{2}=\frac{1}{\left(  c_{A}^{\prime}\left(  \lambda_{k}\right)
\right)  ^{2}}\sum_{j=1}^{N}\det\left(  A-\lambda_{k}I\right)
_{\operatorname{row}j=b}^{2} \label{beta_k^2_sum_squares}%
\end{equation}
Using $\sum_{j=1}^{N}\left(  x_{k}\right)  _{j}^{2}=1$ in combination with
(\ref{x_k_j_detA/j}) yields%
\begin{equation}
\frac{1}{\beta_{k}^{2}}=\sum_{j=1}^{N}\frac{\det^{2}\left(  A_{G\backslash
\left\{  j\right\}  }-\lambda_{k}I\right)  }{\det^{2}\left(  A-\lambda
_{k}I\right)  _{\operatorname{row}j=b}} \label{1_over_beta_k^2}%
\end{equation}
Finally, it follows from the Cauchy-Schwarz inequality applied to
(\ref{beta_k^2}) (or to (\ref{Sum=1})) and with (\ref{beta_k^2_sum_squares})
(or with (\ref{1_over_beta_k^2})) that $\beta_{k}^{2}\leq\left\Vert
b\right\Vert _{2}^{2}$, which leads to the same bound as the 2-norm of a
vector $\left\Vert y\right\Vert _{2}^{2}=y^{T}y$, namely $\beta_{k}%
^{2}=\left(  b^{T}x_{k}\right)  ^{2}\leq\left\Vert b\right\Vert _{2}%
^{2}\left\Vert x_{k}\right\Vert _{2}^{2}=\left\Vert b\right\Vert _{2}^{2}$.

\end{document}